\documentclass[aap,MSNbibl,seceqn,secthm,citesort,dvips]{arximspdf}

\doi{10.1214/10-AAP688}
\volume{20}
\issue{6}
\pubyear{2010}
\firstpage{2318}
\lastpage{2345}

\makeatletter

\newtheorem{lem}[thm]{Lemma}
\newproclaim{defn}[thm]{Definition}
\newproclaim{condkr}{Condition}
\newproclaim{condc}{Condition}
\newproclaim{condn}{Condition}
\newproclaim{condo}{Condition}
\newproclaim{conduo}{Condition}
\newproclaim{condk}{Condition}
\newproclaim{condbb}{Condition}
\newproclaim{condb}{Condition}
\newproclaim{rem}[thm]{Remark}
\newproclaim{example}[thm]{Example}

\newcommand{\pimax}{\pi^{\max}}
\newcommand{\pimin}{\pi^{\min}}
\newcommand{\xrightarrow}[1]{\stackrel{#1}{\longrightarrow}}
\newcommand{\supp}{\operatorname{supp}}

\makeatother

\begin{document}
\begin{frontmatter}

\title{A complete solution to Blackwell's unique ergodicity problem for
hidden Markov chains}
\runtitle{Unique ergodicity and hidden Markov chains}

\begin{aug}
\author[A]{\fnms{Pavel} \snm{Chigansky}\ead[label=e1]{pchiga@mscc.huji.ac.il}\thanksref{aut1}}
\and
\author[B]{\fnms{Ramon} \snm{van Handel}\ead[label=e2]{rvan@princeton.edu}\corref{}}

\thankstext{aut1}{Supported by ISF Grant 314/09.}

\runauthor{P.\ Chigansky and R.\ van Handel}

\affiliation{Hebrew University and Princeton University}

\address[A]{Department of Statistics\\
Hebrew University\\
Mount Scopus\\
Jerusalem 91905\\
Israel \\ \printead{e1}}

\address[B]{Department of Operations Research\\
\quad and Financial Engineering\\
Princeton University\\
Princeton, New Jersey 08544 \\
USA \\ \printead{e2}}
\end{aug}

\received{\smonth{10} \syear{2009}}

%
\begin{abstract}
We develop necessary and sufficient conditions for uniqueness of
the invariant measure of the filtering process associated to an ergodic
hidden Markov model in a finite or countable state space. These results
provide a complete solution to a problem posed by Blackwell (1957), and
subsume earlier partial results due to Kaijser, Kochman and Reeds. The
proofs of our main results are based on the stability theory of nonlinear
filters.
\end{abstract}



\end{frontmatter}

\section{Introduction}
\label{sec:intro}

The interest in the stationary behavior of hidden Markov models dates back
at least to a 1957 paper by Blackwell \cite{Bla57}, who was motivated by
the following problem from information theory. Suppose that $(X_n)_{n\ge
0}$ is a stationary Markov chain which takes values in a finite set. The
entropy rate of such a chain admits a simple expression in terms of its
transition probabilities and stationary distribution. The purpose of the
paper by Blackwell was to obtain a similar expression for the entropy rate
of the stochastic process $Y_n=h(X_n)$, where $h$ is a noninvertible
function. The latter expression does not involve directly the stationary
distribution of the process $(X_n)_{n\ge0}$, but rather a particular
stationary distribution of the associated filtering process $(\pi
_n)_{n\ge
0}$, which is a measure-valued Markov process defined as
$\pi_n:=\mathbf{P}(X_n\in\cdot|Y_1,\ldots,Y_n)$.

The result of Blackwell raises a natural question: does the filtering
process possess a unique stationary measure or, in other words, is the
filtering process \emph{uniquely ergodic}? Blackwell conjectured that the
filter is uniquely ergodic, provided that the underlying Markov chain
$(X_n)_{n\ge0}$ is irreducible. However, as is pointed out by Kaijser
\cite{Kai75}, one of Blackwell's own counterexamples demonstrates that
this conjecture is incorrect. The problem of finding a complete
characterization of the unique ergodicity of the filtering process has
hitherto remained open. The present paper provides one solution to this
problem (in a more general setting).

\subsection{The contributions of Kaijser, Kochman and Reeds}

To our knowledge, the only direct contributions to the problem studied in
this paper are contained in Blackwell's 1957 paper \cite{Bla57}, in a 1975
paper by Kaijser \cite{Kai75} and in two recent papers by Kochman and
Reeds \cite{KR06} and by Kaijser \cite{Kai09}, which we presently review.

In the 1975 paper \cite{Kai75}, Kaijser observes that the filtering
process can be expressed as the ratio of two quantities which are defined
in terms of products of random matrices. Therefore, the unique ergodicity
problem can be studied by means of the Furstenberg--Kesten theory of random
matrix products. Such an analysis leads Kaijser to introduce a certain
subrectangularity condition on the matrices that define the filter
[Condition \ref{condK} in Section \ref{sec:kaijser}]. This rather strong
condition is sufficient, but not necessary for unique ergodicity.
It should be noted that Blackwell's original paper \cite{Bla57} already
gives a sufficient condition for unique ergodicity, which is, however, even
stronger than Kaijser's subrectangularity condition.

In their 2006 paper \cite{KR06}, Kochman and Reeds introduce a weaker
sufficient condition for unique ergodicity of the filter, which requires
that the closure of a certain cone of matrices contains an element of rank
one [Condition \ref{condKR} in Section \ref{sec:mr}]. Kochman and Reeds
demonstrate by means of an explicit computation that Kaijser's condition
implies the rank one condition, but a counterexample shows that the latter
condition is strictly weaker. Besides providing a generalization of
Kaijser's result, Kochman and Reeds employ a different method of proof
that is based on a general result in the ergodic theory of Markov chains
in topological state spaces (which is applied to the filtering process).

Finally, in a recent paper \cite{Kai09}, Kaijser presents an extension of
the result of Kochman and Reeds to hidden Markov models where the
underlying Markov chain $(X_n)_{n\ge0}$ takes values in a countable state
space. (It should be noted that Kochman and Reeds, as well as Kaijser,
admit a more general observation structure than in Blackwell's original
problem.) The extension is far from straightforward, as the ergodic theory
employed by Kochman and Reeds is restricted to Markov chains in locally
compact state spaces, while the space of probability measures on a
countable set is certainly not locally compact. A large part of this
lengthy paper is taken up with the development of a rather specialized
ergodic theorem for Markov chains in Polish spaces, from which a condition
similar in spirit to Kochman and Reeds' rank one condition [Condition
\ref{condB1} in Section \ref{sec:kaijser}] can be derived.

\subsection{The approach of Kunita and filter stability}

Independently from Blackwell's unique ergodicity problem, a general study
of the ergodic theory of nonlinear filtering processes was initiated in
the seminal 1971 paper of Kunita \cite{Kun71}. Kunita studies a somewhat
different problem, in continuous time and with white noise type
observations, but which otherwise bears strong similarities to the problem
studied by Blackwell. In contrast to the approaches developed by
Kaijser, Kochman and Reeds, who study the \emph{equations} that define the
filter using general methods (products of random matrices and ergodic
theory of Markov chains), Kunita studies the nonlinear filter directly
through its characterization as a conditional expectation (an approach we
called \emph{intrinsic} in \cite{CV10}). The techniques developed by
Kunita are in fact extremely general and can be applied also to
Blackwell's problem, though this approach has not previously been
systematically exploited.

Kunita characterizes the invariant measures of the filtering process by
means of the convex ordering. When the signal $(X_n)_{n\ge0}$ is
uniquely ergodic, all invariant measures of the filter are sandwiched
between two distinguished invariant measures which are minimal and maximal
with respect to the convex order, respectively (see Remark
\ref{rem:explkunita} below for a more precise statement). The filter is
uniquely ergodic precisely when the minimal and maximal invariant measures
coincide. The main result of Kunita's paper claims that this is always
the case, when the signal is ergodic in a certain sense. Unfortunately,
the proof of this result contains a serious gap \cite{BCL04}; indeed, the
correctness of the proof is already contradicted by the counterexample
given in Kaijser \cite{Kai75} (see \cite{BCL04,Bud03} for extensive
discussion).

The gap in Kunita's main result is now largely resolved \cite{Van08}, but
under an additional \emph{nondegeneracy} assumption on the observation
structure [Condition~\ref{condN} of Section \ref{sec:kaijser} in the present
setting]. This assumption holds, for example, if
$Y_n=h(X_n)+\varepsilon\xi_n$ where $\xi_n$ is an independent Gaussian
random variable and $\varepsilon>0$ is an arbitrarily small noise
strength, but breaks down in the noiseless case $\varepsilon=0$. The
nondegeneracy assumption evidently captures the phenomenon that
observation noise has a stabilizing effect on the filter, as is the
case in a large number of interesting applications. Unfortunately, it is
the degenerate case that is chiefly of interest in Blackwell's problem,
and unique ergodicity turns out to be more delicate in this setting as
is demonstrated by various counterexamples \cite{Kai75,BCL04,KR06}.

In recent years, there has been considerable interest in the somewhat
different problem of \emph{filter stability} (see the survey \cite{CV10}).
Roughly speaking, the filtering process is called stable if $\pi_n$
becomes independent of its initial condition $\pi_0$ as $n\to\infty$ in a
certain pathwise sense (e.g., as in Theorem \ref{thm:kunita}).
However, it is now well established that when the signal $(X_n)_{n\ge0}$
is ergodic, filter stability and unique ergodicity of the filter are
essentially equivalent properties \cite{Bud03,DS05,Ste89}. In the present
setting, this has two important consequences. First, filter stability can
be used as a tool to study unique ergodicity of the filter, a fact that is
heavily exploited in this paper. Second, previous work on the filter
stability problem provides a set of sufficient conditions for Blackwell's
unique ergodicity problem which are distinct from those proposed by
Kaijser, Kochman and Reeds.

\subsection{Contributions of this paper}

The present paper was inspired by the observation that the conditions of
Kochman and Reeds \cite{KR06} and Kaijser \cite{Kai09} are reminiscent of
the filter stability property, albeit along a single sample path. It is
therefore a natural step to make the connection with filter stability
theory and Kunita's ergodic theory. Our results demonstrate that this
approach is both natural and fruitful.

Our first main result, Theorem \ref{thm:mainc}, establishes that a certain
Condition \ref{condC} is necessary and sufficient for unique ergodicity of the
filter in the case where $X_n$ and $Y_n$ both take values in an (at most)
countable state space. It is easily shown, as we do in Section
\ref{sec:kaijser}, that the sufficient conditions given in Kaijser's
recent paper \cite{Kai09} imply Condition \ref{condC}. It should be noted that
the proof of Theorem \ref{thm:mainc} is surprisingly easy and
natural---that is, provided the connection between filter stability and
Kunita's ergodic theory (given in Theorem \ref{thm:kunita}) is taken for
granted.

Our second main result, Theorem \ref{thm:mainf}, shows that the rank one
Condition \ref{condKR} of Kochman and Reeds is necessary and sufficient for unique
ergodicity of the filter in the case where $X_n$ takes values in a finite
state space. Sufficiency was already proved by Kochman and Reeds
\cite{KR06}, though we give here an entirely different proof of this fact
by showing that Condition \ref{condKR} implies Condition \ref{condC}. The necessity of
Condition \ref{condKR} is new, and answers in the affirmative the question posed
on the last page of \cite{KR06}. Thus the necessity and sufficiency
of Condition~\ref{condKR} provides a complete solution to the original problem
posed by Blackwell \cite{Bla57}.

Our main results subsume all of the sufficient conditions introduced in
the papers of Kaijser, Kochman and Reeds. In addition, we discuss in
Section \ref{sec:kaijser} some sufficient conditions of a different nature
which are inherited from results in the filter stability literature.
Though these conditions are not necessary, they may be substantially
easier to check in practice than Condition \ref{condC} or \ref{condKR}. Moreover, such
conditions remain of independent interest, as we were not able to verify
by an explicit computation that they imply Condition \ref{condC} or \ref{condKR} (of
course, this implication follows indirectly from the necessity of these
conditions).

\subsection{Organization of the paper}

The remainder of this paper is organized as follows. In Section
\ref{sec:prelim} we introduce the basic hidden Markov model, and we fix
once and for all the notation and standing assumptions that are presumed
to be in force throughout the paper. We also state our main results,
Theorems \ref{thm:mainc} and~\ref{thm:mainf}. In Section
\ref{sec:ergodic}, we introduce the connection between filter stability
and unique ergodicity of the filter. The main result of this section,
Theorem \ref{thm:kunita}, adapts the necessary theory to the setting of
this paper and forms the foundation for the proofs of our main results.
Section \ref{sec:mainc} is devoted to the proof of Theorem
\ref{thm:mainc}, while Section \ref{sec:mainf} is devoted to the proof of
Theorem \ref{thm:mainf}. Section \ref{sec:kaijser} develops various
sufficient conditions for unique ergodicity within the setting of this
paper. Finally, the \hyperref[app:proofs]{Appendix} is devoted to the proofs of various results
that were omitted from the body of the paper.

\section{Preliminaries and main results}
\label{sec:prelim}

\subsection{The canonical setup and standing assumptions}

Throughout this paper, we operate in the following setup. We consider
the stochastic process $(X_n,Y_n)_{n\in\mathbb{Z}}$, where $X_n$
takes values in the state space $E$, and $Y_n$ takes values in the state
space~$F$. We will always presume that the following assumptions are in
force:
\begin{itemize}
\item$E$ is either finite ($E=\{1,\ldots,p\}$) or countable
($E=\mathbb{N}$).
\item$F$ is a Polish space [endowed with its Borel $\sigma$-field
$\mathcal{B}(F)$].
\end{itemize}
We realize the stochastic process $(X_n,Y_n)_{n\in\mathbb{Z}}$ on the
canonical path space $\Omega=\Omega^X\times\Omega^Y$ with
$\Omega^X=E^{\mathbb{Z}}$ and $\Omega^Y=F^{\mathbb{Z}}$, such that
$X_n(x,y)=x(n)$ and $Y_n(x,y)=y(n)$. Denote by $\mathcal{F}$ the Borel
$\sigma$-field on $\Omega$, and introduce the $\sigma$-fields
\[
\mathcal{F}_{m,n}^X = \sigma\{X_k\dvtx k\in[m,n]\},\qquad
\mathcal{F}_{m,n}^Y = \sigma\{Y_k\dvtx k\in[m,n]\}
\]
for $m,n\in\mathbb{Z}$, $m\le n$. The $\sigma$-fields
$\mathcal{F}_{-\infty,n}^X$, $\mathcal{F}_{m,\infty}^X$, etc., are defined
in the usual fashion (e.g., $\mathcal{F}_{-\infty,n}^X=
\bigvee_{m\le n}\mathcal{F}_{m,n}^X$). For future reference, we define
\[
\mathcal{G}_{m,n} = \mathcal{F}^X_{-\infty,m}\vee
\mathcal{F}^Y_{-\infty,n},\qquad
\mathcal{G}_{-\infty,n} = \bigcap_{m\le n}\mathcal{G}_{m,n}
\]
(note that $\mathcal{F}^Y_{-\infty,n}\subset\mathcal{G}_{-\infty,n}$,
a fact that will be used frequently in the following).
Finally, the shift $\Theta\dvtx\Omega\to\Omega$ is defined
as $\Theta(x,y)(m) = (x(m+1),y(m+1))$.

We now define a probability measure on $(\Omega,\mathcal{F})$ under which
$(X_k,Y_k)_{k\in\mathbb{Z}}$ is a \emph{hidden Markov model}. Our model
is specified by the following ingredients:
\begin{enumerate}[(2)]
\item[(1)] A $\sigma$-finite reference measure $\varphi$ on $F$.
\item[(2)] A nonnegative matrix function $M\dvtx F\to\mathbb{R}_+^{E\times E}$
such that
\[
\sup_{i\in E}\sum_{j\in E}M_{ij}(y)<\infty\qquad
\mbox{for }\varphi\mbox{-a.e.\ }y\in F,
\]
and such that the matrix
\[
P=(P_{ij})_{i,j\in E},\qquad
P_{ij}:= \int M_{ij}(y) \varphi(dy)
\]
defines the transition matrix of an irreducible and positive recurrent
(but not necessarily aperiodic) Markov chain in the state space $E$.
\end{enumerate}

As $P$ is irreducible and positive recurrent, there is a unique
probability measure $\lambda$ on $E$ that is invariant $\lambda P =
\lambda$ (as is usual, we identify measures and functions on a countable
space with row and column vectors, respectively). A standard extension
argument allows us to construct a probability measure $\mathbf{P}$ on
$(\Omega,\mathcal{F})$ under which $(X_k,Y_k)_{k\in\mathbb{Z}}$ is a
stationary Markov chain with transition probabilities
\[
\mathbf{P}(X_k=j, Y_k\in A|X_{k-1}=i, Y_{k-1}=y) =
\int_A M_{ij}(y') \varphi(dy')
\]
for $i,j\in E$, $y\in F$, $A\in\mathcal{B}(F)$. It should be noted that
under $\mathbf{P}$, the process $(X_k)_{k\in\mathbb{Z}}$ is a stationary
Markov chain with transition matrix $P$, and $(Y_k)_{k\in\mathbb{Z}}$ are
conditionally independent given $(X_k)_{k\in\mathbb{Z}}$. This is
precisely the defining property of a hidden Markov model. The process
$(X_k)_{k\in\mathbb{Z}}$ represents an unobserved or ``hidden'' signal
process, while $(Y_k)_{k\in\mathbb{Z}}$ is the observation process.
The canonical probability space $(\Omega,\mathcal{F},\mathbf{P})$
thus constructed will remain fixed throughout the paper.

\begin{rem}
A hidden Markov model is often assumed to satisfy the additional
assumption that $Y_k$ is a (noisy) function of $X_k$ only. In this case,
one can factor $M_{ij}(y)=P_{ij}R_j(y)$, where $R_j(y)$ is the density of
$\mathbf{P}(Y_k\in\cdot|X_k=j)$ with respect to $\varphi$. In the
present setting, the conditional law of $Y_k$ can depend on both $X_k$
and $X_{k-1}$. The generalization afforded by this model is minor, but
allows us to include the partitioned transition matrices of
\cite{KR06,Kai09} as a special case.
\end{rem}

\begin{rem}
The boundedness condition $\sup_{i\in E}\sum_{j\in E}M_{ij}(y)<\infty$
a.e.\ is automatically satisfied in the following cases:
\begin{itemize}
\item When $E$ is a finite set, the condition holds trivially.
\item When $E$ is countable and $F$ is at most countable, the condition
always holds. Indeed, note that in this case $\sum_{y\in F}\sum_{j\in E}
M_{ij}(y) \varphi(\{y\})= \sum_{j\in E} P_{ij}=1$, so that
$\sup_{i\in E}\sum_{j\in E}M_{ij}(y)\le\varphi(\{y\})^{-1}
<\infty$ for $\varphi$-a.e.\ $y\in F$.
\end{itemize}
The significance of this assumption is that it ensures the Feller property
of the filter.
\end{rem}

For any Polish space $S$ we denote by $\mathcal{B}(S)$ the Borel
$\sigma$-field of $S$, by $\mathcal{P}(S)$ the space of probability
measures on $S$, and by $\mathcal{C}_b(S)$ the space of bounded continuous
functions on $S$. We will always endow $\mathcal{P}(S)$ with the topology
of weak convergence of probability measures [recall that $\mathcal{P}(S)$
is then itself Polish], and we write $\mu_n\Rightarrow\mu$ if the sequence
$(\mu_n)\subset\mathcal{P}(S)$ converges weakly to $\mu\in\mathcal{P}(S)$.
The total variation distance between probability measures
$\mu,\nu\in\mathcal{P}(S)$ is defined as
\[
\|\mu-\nu\| = \sup_{\|f\|_\infty\le1}
\biggl|\int f \,d\mu- \int f \,d\nu\biggr|.
\]
Finally, let us recall that as $E$ is at most countable and $P$ is
irreducible, the invariant measure $\lambda$ must charge every point of
$E$. Therefore $\mu\ll\lambda$ for every $\mu\in\mathcal{P}(E)$, and we
can define the probability measures $\mathbf{P}^\mu$ on
$(\Omega,\mathcal{F})$ as
\[
\frac{d\mathbf{P}^\mu}{d\mathbf{P}} =
\frac{d\mu}{d\lambda}(X_0),\qquad
\mu\in\mathcal{P}(E).
\]
The restriction of $\mathbf{P}^\mu$ to $\mathcal{F}_{0,\infty}^X
\vee\mathcal{F}_{1,\infty}^Y$ defines a hidden Markov model with the same
transition probabilities as under $\mathbf{P}$, but with the
initial distribution $X_0\sim\mu$. If the initial distribution is a point
mass on $x\in E$, we will write $\mathbf{P}^x$ instead of
$\mathbf{P}^{\delta_x}$.

\subsection{Nonlinear filtering}

The purpose of nonlinear filtering is to compute the conditional
distribution of the hidden signal given the available observations.
In this paper we will encounter several variants of the nonlinear filter,
defined as follows:
\[
\pi_n = \mathbf{P}(X_n\in\cdot|\mathcal{F}_{1,n}^Y),\qquad
\pi_n^\mu= \mathbf{P}^\mu(X_n\in\cdot|\mathcal{F}_{1,n}^Y),\qquad
\pi_n^x = \mathbf{P}^x(X_n\in\cdot|\mathcal{F}_{1,n}^Y)
\]
for $n\in\mathbb{Z}_+$, $\mu\in\mathcal{P}(E)$, $x\in E$ (here
$\pi_0=\lambda$, $\pi_0^\mu=\mu$ and $\pi_0^x=\delta_x$) and
\[
\pimin_n = \mathbf{P}(X_n\in\cdot|\mathcal{F}_{-\infty,n}^Y),\qquad
\pimax_n = \mathbf{P}(X_n\in\cdot|\mathcal{G}_{-\infty,n})
\]
for $n\in\mathbb{Z}$. Though the relevance of $\pimin_n$ and
$\pimax_n$ may not be entirely evident at present, their role will be
clarified in Section \ref{sec:ergodic} below.

The following elementary results are essentially
known; short proofs are provided in Appendix \ref{sec:basicfiltering} for
the reader's convenience.

\begin{lem}[(Filtering recursion)]
\label{lem:filtrecurs}
For any $m,n\in\mathbb{Z}$, $n>m$ we have $\mathbf{P}$-a.s.
\[
\pimin_n = \frac{\pimin_m M(Y_{m+1})\cdots M(Y_n)}{
\pimin_m M(Y_{m+1})\cdots M(Y_n)1},\qquad
\pimax_n = \frac{\pimax_m M(Y_{m+1})\cdots M(Y_n)}{
\pimax_m M(Y_{m+1})\cdots M(Y_n)1}.
\]
Similarly, for any $n>m\ge0$, we have $\mathbf{P}^\mu$-a.s.
\[
\pi_n^\mu= \frac{\pi_m^\mu M(Y_{m+1})\cdots M(Y_n)}{
\pi_m^\mu M(Y_{m+1})\cdots M(Y_n)1}.
\]
The recursion for $\pi_n,\pi_n^x$ is obtained by
choosing $\mu=\lambda$ or $\mu=\delta_x$, respectively.
\end{lem}

It should be noted that $\pi_n^\mu$ is defined only up to a
$\mathbf{P}^\mu$-null set. Indeed,
\[
\mathbf{P}^\mu\bigl((Y_1,\ldots,Y_n)\in A\bigr) =
\int_A \mu M(y_1)\cdots M(y_n)1 \varphi(dy_1)\cdots\varphi(dy_n),
\]
that is, $\mu M(y_1)\cdots M(y_n)1$ is the density of the law of
$(Y_1,\ldots,Y_n)$ under $\mathbf{P}^\mu$. Similarly, it is easily seen
that $\pi_m^\mu M(y_{m+1})\cdots M(y_n)1$ is the density of the law of
$(Y_{m+1},\ldots,Y_n)$ under the conditional measure
$\mathbf{P}^\mu( \cdot|Y_0,\ldots,Y_m)$. Therefore, the denominator in
the filtering recursion can only vanish on a $\mathbf{P}^\mu$-null set.
Similar considerations hold for $\pimin_n,\pimax_n$, which are defined up
to a $\mathbf{P}$-null set.

\begin{lem}[(Markov property)]
\label{lem:filtmarkov}
$(\pimin_n)_{n\in\mathbb{Z}}$, $(\pimax_n)_{n\in\mathbb{Z}}$
are stationary $\mathcal{P}(E)$-valued Markov chains under $\mathbf{P}$,
whose transition kernel $\mathsf{\Pi}$ is defined by
\[
\int f(\nu) \mathsf{\Pi}(\mu,d\nu) =
\int f \biggl(\frac{\mu M(y)}{\mu M(y)1}
\biggr) \mu M(y)1 \varphi(dy),\qquad
f\in\mathcal{C}_b(\mathcal{P}(E)).
\]
Similarly, $(\pi_n^\mu)_{n\in\mathbb{Z}_+}$ is a Markov chain under
$\mathbf{P}^\mu$ with transition kernel $\mathsf{\Pi}$.
\end{lem}

\begin{rem}
As $(\pimin_n)_{n\in\mathbb{Z}}$, $(\pimax_n)_{n\in\mathbb{Z}}$ are
stationary Markov chains with transition kernel $\mathsf{\Pi}$, the laws
of $\pimax_0$ and $\pimin_0$ must be invariant for $\mathsf{\Pi}$.
Therefore, the filter always possesses at least one invariant measure.
\end{rem}

\subsection{Main results}
\label{sec:mr}

This paper aims to resolve the following question: when does the filter
possess a \emph{unique} invariant measure, that is, when does the equation
$\mathsf{M}\mathsf{\Pi}=\mathsf{M}$ possess a unique solution
$\mathsf{M}\in\mathcal{P}(\mathcal{P}(E))$?

We begin by establishing a general sufficient condition for unique
ergodicity, which is also necessary when the observation state space $F$
is at most countable.

\renewcommand{\thecondc}{(C)}
\begin{condc}\label{condC}
For every $\varepsilon>0$, there exist an integer $N\in\mathbb{N}$
and subsets $\mathcal{S}\subset\mathcal{P}(E)$ and $\mathcal{O}\subset
F^N$ such that the following hold:
\begin{enumerate}[(3)]
\item[(1)]$\mathbf{P}(\pimin_0\in\mathcal{S}\ \mathrm{and}\
\pimax_0\in\mathcal{S})>0$ and $\varphi^{\otimes N}(\mathcal{O})>0$.
\item[(2)]$\mu M(y_1)\cdots M(y_N)1>0$ for all
$\mu\in\mathcal{S}$ and $(y_1,\ldots,y_N)\in\mathcal{O}$.
\item[(3)] For all $\mu,\nu\in\mathcal{S}$ and $(y_1,\ldots,y_N)\in\mathcal{O}$
\[
\biggl\|\frac{\mu M(y_1)\cdots M(y_N)}{\mu M(y_1)\cdots M(y_N)1}
-\frac{\nu M(y_1)\cdots M(y_N)}{\nu M(y_1)\cdots M(y_N)1}
\biggr\|<\varepsilon.
\]
\end{enumerate}
\end{condc}

\begin{thm}
\label{thm:mainc}
Suppose that Condition \ref{condC} holds. Then the filter admits a unique
invariant measure $\mathsf{M}$, and we have
$n^{-1}\sum_{k=1}^n\mathsf{M}_0\mathsf{\Pi}^k\Rightarrow\mathsf{M}$
as $n\to\infty$ for any $\mathsf{M}_0\in\mathcal{P}(\mathcal{P}(E))$.
If, in addition, the signal transition matrix $P$ is aperiodic, then
we have $\mathsf{M}_0\mathsf{\Pi}^n\Rightarrow\mathsf{M}$ as $n\to\infty$
for any $\mathsf{M}_0\in\mathcal{P}(\mathcal{P}(E))$.

Conversely, suppose that the observation state space $F$ is a finite or
countable set, and that the filter is uniquely ergodic. Then
Condition \ref{condC} holds.
\end{thm}

The proof of this result is given in Section \ref{sec:mainc}.

Next, we consider the following condition, due to Kochman and Reeds
\cite{KR06}, for the case where the signal state space $E$ is a finite
set.

\renewcommand{\thecondkr}{(KR)}
\begin{condkr}\label{condKR}
Let $E$ be a finite set, and define the cone of matrices
\[
\mathcal{K} = \{c M(y_1)\cdots M(y_n)\dvtx n\in\mathbb{N},
y_1,\ldots,y_n\in F, c\in\mathbb{R}_+\}.
\]
Then the closure $\operatorname{cl} \mathcal{K}$ contains a matrix of rank 1.
\end{condkr}

Kochman and Reeds prove that this condition is sufficient for uniqueness
of the invariant measure of the filter (in \cite{KR06}, both $E$ and $F$
are presumed to be finite). The following result shows that Condition \ref{condKR} is in fact \emph{equivalent} to unique ergodicity of the filter, as
well as to Condition \ref{condC} above, when the signal state space is a finite
set. This provides a complete solution to a problem posed by Blackwell
\cite{Bla57}, and answers in the affirmative the question posed at the end
of \cite{KR06}.

\begin{thm}
\label{thm:mainf}
Suppose $E$ is a finite set and that one of the following hold:
\begin{itemize}
\item$F$ is a finite or countable set, and $\varphi$ is the counting
measure; or
\item$F=\mathbb{R}^d$, $\varphi$ is the Lebesgue measure,
and $y\mapsto M(y)$ is continuous.
\end{itemize}
Then the following are equivalent:
\begin{enumerate}[(3)]
\item[(1)] The filter admits a unique invariant measure $\mathsf{M}$.
\item[(2)] Condition \ref{condKR} holds.
\item[(3)] Condition \ref{condC} holds.
\end{enumerate}
When any of these conditions hold, we have $n^{-1}\sum_{k=1}^n\mathsf{M}_0
\mathsf{\Pi}^k\Rightarrow\mathsf{M}$ as $n\to\infty$ for any
$\mathsf{M}_0\in\mathcal{P}(\mathcal{P}(E))$. If, in addition, the signal
transition matrix $P$ is aperiodic, then we have
$\mathsf{M}_0\mathsf{\Pi}^n\Rightarrow\mathsf{M}$ as $n\to\infty$
for any $\mathsf{M}_0\in\mathcal{P}(\mathcal{P}(E))$.
\end{thm}

The proof will be given in Section \ref{sec:mainf}.

Finally, various sufficient conditions for unique ergodicity of the filter
were given by Kaijser \cite{Kai75,Kai09}. These conditions are easily
shown to imply Condition~\ref{condC}, as is discussed in Section
\ref{sec:kaijser}. We therefore reproduce Kaijser's results using a much
simpler proof. Similarly, various conditions that have been introduced in
the context of filter stability \cite{Van08,Van09,CV10} are shown in
Section \ref{sec:kaijser} to imply unique ergodicity of the filter. None
of the latter sufficient conditions is also necessary; however, when they
apply, they are often easier to check than Condition \ref{condC} or \ref{condKR}.

\section{Ergodic theory and stability of nonlinear filters}
\label{sec:ergodic}

The proofs of our main results are based on a general circle of ideas
connecting the ergodic theory \cite{Kun71,Ste89} and asymptotic stability
\cite{Van08,CV10} of nonlinear filters. Indeed, it is by now well
established \cite{Bud03,DS05} that unique ergodicity and stability of the
filter are essentially equivalent properties. The purpose of this section
is to introduce the relevant results in this direction that will be needed
in what follows. Though the results in this section are adapted to the
setting of this paper, their proofs largely follow along the lines of
\cite{Kun71,Ste89,DS05,Van08}. We have therefore relegated the proofs to
the \hyperref[app:proofs]{Appendix}.

The following characterization will be of central importance.

\begin{thm}
\label{thm:kunita}
Consider the following conditions:
\begin{enumerate}[(3)]
\item[(1)] The filter possesses a unique invariant measure $\mathsf{M}\in
\mathcal{P}(\mathcal{P}(E))$.
\item[(2)]$\pimax_0=\pimin_0$ $\mathbf{P}$-a.s.
\item[(3)]$\|\pi_n^\mu-\pi_n^\nu\|\to0$ as $n\to\infty$
$\mathbf{P}^\mu$-a.s.\ whenever $\mu\ll\nu$.
\item[(4)]$n^{-1}\sum_{k=1}^n\mathsf{M}_0
\mathsf{\Pi}^k\Rightarrow\mathsf{M}$ as $n\to\infty$ for any
$\mathsf{M}_0\in\mathcal{P}(\mathcal{P}(E))$.
\item[(5)]$\mathsf{M}_0\mathsf{\Pi}^n\Rightarrow\mathsf{M}$ as $n\to\infty$
for any $\mathsf{M}_0\in\mathcal{P}(\mathcal{P}(E))$.
\end{enumerate}
Conditions \textup{1--4} are equivalent. If, in addition, the signal
transition matrix $P$ is aperiodic, then conditions \textup{1--5} are equivalent.
\end{thm}

The proof is given in Appendix \ref{sec:kunita}.

\begin{rem}
\label{rem:explkunita}
Condition 1 is the desired unique ergodicity property of the filter.
Condition 3 is the filter stability property. Conditions 4 and 5
characterize the convergence of the law of the filter to the invariant
measure.

Condition 2 in Theorem \ref{thm:kunita} stems from an ingenious device
introduced by Kunita in the seminal paper \cite{Kun71} and used in the
proof of Theorem \ref{thm:kunita}. By Lemma~\ref{lem:filtmarkov},
$(\pimin_n)_{n\in\mathbb{Z}}$ and $(\pimax_n)_{n\in\mathbb{Z}}$
are stationary Markov processes. Therefore, the laws $\mathsf{M}^{\max},\mathsf{M}^{\min}\in\mathcal{P}(\mathcal{P}(E))$ of the
$\mathcal{P}(E)$-valued random variables $\pimax_0,\pimin_0$ are
invariant for the filter transition kernel $\mathsf{\Pi}$.
Kunita shows that \emph{any} invariant measure $\mathsf{M}$ for
$\mathsf{\Pi}$ is sandwiched between $\mathsf{M}^{\max}$ and
$\mathsf{M}^{\min}$ in the sense that
\[
\int f(\mu) \mathsf{M}^{\min}(d\mu) \le
\int f(\mu) \mathsf{M}(d\mu) \le
\int f(\mu) \mathsf{M}^{\max}(d\mu)
\]
for every convex function $f\in\mathcal{C}_b(\mathcal{P}(E))$. In other
words, within the family of $\mathsf{\Pi}$-invariant measures,
$\mathsf{M}^{\min}$ is minimal and $\mathsf{M}^{\max}$ is maximal
with respect to the convex ordering. The identity $\pimax_0=\pimin_0$
ensures that the maximal and minimal invariant measures are identical, so
that there can be only one invariant measure.
\end{rem}

\begin{example}
Some intuition may be obtained from the following simple example
\cite{KR06}, which is a typical case where the filter fails to
be uniquely ergodic. Let $E=F=\{0,1\}$ (endowed with the counting
measure), and let
\[
M(0) =
\pmatrix{
0 & 1/2 \cr
1/2 & 0
}
,\qquad
M(1) =
\pmatrix{
1/2 & 0 \cr
0 & 1/2
}
.
\]
Note that $P=M(0)+M(1)$ is irreducible and aperiodic with invariant
measure $\lambda_0=\lambda_1=1/2$, and $Y_k=I_{\{X_{k-1}=X_k\}}$ for
all $k\ge1$.

As $(Y_k)_{k\ge1}$ reveals exactly when the transitions of
$(X_k)_{k\ge0}$ occur, we evidently have $X_k\in\sigma\{X_{m},Y_{m+1},
Y_{m+2},\ldots,Y_n\}$ for every $m\le k\le n$. It follows that
$\mathcal{G}_{m,n} = \mathcal{F}^X_{-\infty,n}$ for every $m\le n$, so
that in particular $\mathcal{G}_{-\infty,n}=\mathcal{F}^X_{-\infty,n}$.
Therefore
\[
\pimax_n = \delta_{X_n},\qquad
\mathsf{M}^{\max} =
\tfrac{1}{2} \{
\delta_{\delta_{0}}+\delta_{\delta_{1}}
\}.
\]
On the other hand, it follows immediately from the filtering recursion
that $\pi_n=\lambda$ for all $n$. It is therefore not difficult to show
that $\pimin_n=\lambda$ also, so that
\[
\pimin_n=\lambda,\qquad
\mathsf{M}^{\min} = \delta_{\lambda}
= \delta_{(\delta_0+\delta_1)/2}.
\]
With a little more work, one can show that any invariant measure
$\mathsf{M}$ is of the form
\[
\mathsf{M} =
\int_0^{1/2}
\frac{\delta_{\epsilon\delta_0+(1-\epsilon)\delta_1}+
\delta_{(1-\epsilon)\delta_0+\epsilon\delta_1}}{2}
m(d\epsilon)
\]
for some probability measure $m$ on $[0,1/2]$. It is easily seen that any
such $\mathsf{M}$ does indeed lie between $\mathsf{M}^{\max}$ and
$\mathsf{M}^{\min}$ in the convex ordering.
\end{example}

Besides the characterization of unique ergodicity in Theorem
\ref{thm:kunita}, we will require the following convergence property
which holds regardless of unique ergodicity.

\begin{lem}
\label{lem:limminmax}
$\lim_{n\to\infty}\|\pimax_n-\pimin_n\|$ exists $\mathbf{P}$-a.s.
\end{lem}

The proof of this result is also given in Appendix \ref{sec:kunita}. Its
relevance is due to the following observation. In order to prove
$\pimax_0=\pimin_0$ (hence unique ergodicity by Theorem \ref{thm:kunita}),
it suffices to show that $\lim_{n\to\infty}\|\pimax_n-\pimin_n\|=0$
$\mathbf{P}$-a.s., as $(\pimin_n)_{n\in\mathbb{Z}}$ and
$(\pimax_n)_{n\in\mathbb{Z}}$ are stationary processes. But by virtue of
Lemma \ref{lem:limminmax}, it then suffices to show only that
$\|\pimax_n-\pimin_n\|$ converges to zero along a sequence of stopping
times. The main idea behind the proof of Theorem \ref{thm:mainc} is that
Condition~\ref{condC} allows us to construct explicitly such a sequence stopping
times.\looseness=1\vadjust{\goodbreak}

\section{\texorpdfstring{Proof of Theorem \protect\ref{thm:mainc}}{Proof of Theorem 2.6}}
\label{sec:mainc}

\subsection{\texorpdfstring{Sufficiency of Condition
\protect\ref{condC}}{Sufficiency of Condition (C)}}

We will need the following lemma.

\begin{lem}
\label{lem:ergodic}
The sequence $(X_k,Y_k)_{k\in\mathbb{Z}}$ is ergodic under $\mathbf{P}$.
\end{lem}

\begin{pf}
As the signal transition matrix $P$ is presumed to be irreducible and
positive recurrent, it is easily established that the pair
$(X_k,Y_k)_{k\in\mathbb{Z}}$ is a Markov process that possesses a unique
invariant measure. This measure is therefore trivially an extreme point
of the set of invariant measures, hence ergodic.
\end{pf}

We will also use the following simple result.

\begin{lem}
\label{lem:pospro}
Let the set $\mathcal{S}\subset\mathcal{P}(E)$ and
$\mathcal{O}\subset F^N$ be as in Condition \ref{condC}.
Then we have $\mathbf{P}(\pimin_0\in\mathcal{S},
\pimax_0\in\mathcal{S}, (Y_1,\ldots,Y_N)\in\mathcal{O})>0$.
\end{lem}

\begin{pf}
Let us write for simplicity $Y=(Y_1,\ldots,Y_N)$. As $\pimin_0$ and
$\pimax_0$ are $\mathcal{G}_{-\infty,0}$-measurable by construction,
we have
\begin{eqnarray*}
&&\mathbf{P}(\pimin_0\in\mathcal{S},
\pimax_0\in\mathcal{S}, Y\in\mathcal{O}) \\
&&\qquad =
\mathbf{E}\bigl(I_\mathcal{S}(\pimin_0)
I_\mathcal{S}(\pimax_0)
\mathbf{P}(Y\in\mathcal{O}|\mathcal{G}_{-\infty,0})\bigr)
\\
&&\qquad =
\mathbf{E} \biggl(I_\mathcal{S}(\pimin_0)
I_\mathcal{S}(\pimax_0)
\int_{\mathcal{O}} \pimax_0 M(y_1)\cdots M(y_N)1 \varphi(dy_1)
\cdots\varphi(dy_N) \biggr).
\end{eqnarray*}
It is now easily seen that Condition \ref{condC} implies the result.
\end{pf}

We now proceed with the proof of the sufficiency part of Theorem
\ref{thm:mainc}. Suppose that Condition \ref{condC} holds, and fix an arbitrary
decreasing sequence $\varepsilon_k\searrow0$. Then for every $k$ we can
find $N_k\in\mathbb{N}$, $\mathcal{S}_k\subset\mathcal{P}(E)$, and
$\mathcal{O}_k\subset F^{N_k}$ such that the properties 1--3 of
Condition \ref{condC} are satisfied. Define the events
\[
A_{n,k} = \{\pimin_n\in\mathcal{S}_k,
\pimax_n\in\mathcal{S}_k, (Y_{n+1},\ldots,Y_{n+N_k})\in\mathcal{O}_k\}.
\]
Then, by the stationarity of
$(X_n,Y_n,\pimin_n,\pimax_n)_{n\in\mathbb{Z}}$
(Lemma \ref{lem:shiftminmax}), we have
\[
\lim_{T\to\infty}
\frac{1}{T}\sum_{n=1}^T I_{A_{n,k}} =
\lim_{T\to\infty}
\frac{1}{T}\sum_{n=1}^T \{I_{A_{0,k}}\circ\Theta^n\}
= \mathbf{P}(A_{0,k}) > 0 \qquad \mathbf{P}\mbox{-a.s.},
\]
where we have used Birkhoff's ergodic theorem together with Lemmas
\ref{lem:ergodic} and~\ref{lem:pospro}.
Thus, for any $k$, the event $A_{n,k}$ occurs at a positive rate,
so that certainly
\[
\mathbf{P} \Biggl(
\bigcap_{k=1}^\infty\limsup_{n\to\infty}A_{n,k}
\Biggr)=1.
\]
Now define the stopping times $\tau_0=0$ and
\[
\tau_k=\min\{n>\tau_{k-1}\dvtx
\pimin_{n-N_k}\in\mathcal{S}_k,
\pimax_{n-N_k}\in\mathcal{S}_k,
(Y_{n-N_k+1},\ldots,Y_{n})\in\mathcal{O}_k\}
\]
for any $k\ge1$. It follows directly that
\[
\mathbf{P}(\tau_k<\infty\mbox{ for all }k) \ge
\mathbf{P} \Biggl(
\bigcap_{k=1}^\infty\limsup_{n\to\infty}A_{n,k}
\Biggr)=1.
\]
Moreover, by Condition \ref{condC} and Lemma \ref{lem:filtrecurs}, we have
\begin{eqnarray*}
\|\pimax_{\tau_k}-\pimin_{\tau_k} \| &=&
\biggl\|
\frac{\pimax_{\tau_k-N_k}M(Y_{\tau_k-N_k+1})\cdots
M(Y_{\tau_k})}{\pimax_{\tau_k-N_k}M(Y_{\tau_k-N_k+1})\cdots
M(Y_{\tau_k})1}
\\
&&\hphantom{\biggl\|}{} -
\frac{\pimin_{\tau_k-N_k}M(Y_{\tau_k-N_k+1})\cdots
M(Y_{\tau_k})}{\pimin_{\tau_k-N_k}M(Y_{\tau_k-N_k+1})\cdots
M(Y_{\tau_k})1}
\biggr\| \le\varepsilon_k
\end{eqnarray*}
for all $k\ge1$ $\mathbf{P}$-a.s. Therefore, Lemma \ref{lem:limminmax}
shows that $\|\pimax_n-\pimin_n\|\to0$ as $n\to\infty$ $\mathbf{P}$-a.s.
But using the stationarity of $(\pimin_n,\pimax_n)_{n\in\mathbb{Z}}$
(Lemma \ref{lem:shiftminmax}) and the dominated convergence theorem,
we find that
\[
\mathbf{E}(\|\pimax_0-\pimin_0\|) =
\mathbf{E}(\|\pimax_n-\pimin_n\|)
\xrightarrow{n\to\infty}0,
\]
so that evidently $\pimax_0=\pimin_0$ $\mathbf{P}$-a.s. The sufficiency
part of Theorem \ref{thm:mainc} now follows immediately from Theorem
\ref{thm:kunita}.

\subsection{\texorpdfstring{Necessity of Condition
\protect\ref{condC}}{Necessity of Condition (C)}}

Throughout this subsection, we assume that the observation state space $F$
is finite or countable, and that the filter possesses a unique invariant
measure. We aim to show that Condition \ref{condC} must hold.

Denote by $\mathsf{M}$ the law of $\pimin_0$. Thus $\mathsf{M}$ is
invariant, hence the unique invariant measure of the filter. Fix an
arbitrary state $i\in E$, and note that
\[
\mathbf{E}((\pimin_0)_i)=\lambda_i>0
\quad \mbox{implies}\quad
\mathbf{P}\bigl((\pimin_0)_i\ge\lambda_i/2\bigr)>0.
\]
Therefore, writing
$\mathcal{R}=\{\mu\in\mathcal{P}(E)\dvtx\mu_i\ge\lambda_i/2\}$, we have
$\mathsf{M}(\mathcal{R})>0$. We can thus define the probability measure
$\mathsf{M}_{\mathcal{R}}( \cdot):=\mathsf{M}( \cdot
\mbox{}\cap\mathcal{R})/\mathsf{M}(\mathcal{R})$.

Now note that, by Theorem \ref{thm:kunita}, we have
$\|\pi_n^i-\pi_n^\mu\|\to0$ $\mathbf{P}^i$-a.s.\ for all\vspace*{1pt}
$\mu\in\mathcal{R}$. In particular, the set of points
$((y_k)_{k\in\mathbb{N}},\mu)\in F^{\mathbb{N}}\times\mathcal{P}(E)$
such that
\[
\biggl\|
\frac{\delta_iM(y_1)\cdots M(y_n)}{
\delta_iM(y_1)\cdots M(y_n)1}-
\frac{\mu M(y_1)\cdots M(y_n)}{
\mu M(y_1)\cdots M(y_n)1}
\biggr\|\xrightarrow{n\to\infty}0
\]
has $\mathbf{P}^i\otimes\mathsf{M}_{\mathcal{R}}$-full measure.
It follows that for $\mathbf{P}^i$-a.e.\ path $(y_k)_{k\in\mathbb{N}}$,
the above convergence holds for $\mathsf{M}_{\mathcal{R}}$-a.e.\ $\mu$.
Therefore, as $F$ is at most countable [so the law of $(Y_1,\ldots,Y_n)$
is atomic for all $n<\infty$], we can certainly find a single sequence
$(\tilde y_k)_{k\in\mathbb{N}}$ with $\mathbf{P}^i(Y_1=\tilde
y_1,\ldots,Y_n=\tilde y_n)>0$ for all $n<\infty$ such that
\[
\biggl\|
\frac{\delta_iM(\tilde y_1)\cdots M(\tilde y_n)}{
\delta_iM(\tilde y_1)\cdots M(\tilde y_n)1}-
\frac{\mu M(\tilde y_1)\cdots M(\tilde y_n)}{
\mu M(\tilde y_1)\cdots M(\tilde y_n)1}
\biggr\|\xrightarrow{n\to\infty}0\qquad
\mbox{for }\mathsf{M}_{\mathcal{R}}\mbox{-a.e. } \mu.
\]
By Egorov's theorem, there is a subset $\mathcal{S}\subset\mathcal{R}$
with $\mathcal{M}_{\mathcal{R}}(\mathcal{S})>0$ such that
\[
\sup_{\mu\in\mathcal{S}}
\biggl\|
\frac{\delta_iM(\tilde y_1)\cdots M(\tilde y_n)}{
\delta_iM(\tilde y_1)\cdots M(\tilde y_n)1}-
\frac{\mu M(\tilde y_1)\cdots M(\tilde y_n)}{
\mu M(\tilde y_1)\cdots M(\tilde y_n)1}
\biggr\|\xrightarrow{n\to\infty}0.
\]
We are now in the position to show that Condition \ref{condC} holds true.
Given $\varepsilon>0$, we first choose the integer $N\in\mathbb{N}$
large enough so that
\[
\sup_{\mu\in\mathcal{S}}
\biggl\|
\frac{\delta_iM(\tilde y_1)\cdots M(\tilde y_N)}{
\delta_iM(\tilde y_1)\cdots M(\tilde y_N)1}-
\frac{\mu M(\tilde y_1)\cdots M(\tilde y_N)}{
\mu M(\tilde y_1)\cdots M(\tilde y_N)1}
\biggr\|\le\frac{\varepsilon}{2}.
\]
We let $\mathcal{S}$ be as above and define the singleton
$\mathcal{O}=\{(\tilde y_1,\ldots,\tilde y_N)\}$.
By Theorem~\ref{thm:kunita}, we have $\pimax_0=\pimin_0$ and therefore
\[
\mathbf{P}(\pimax_0\in\mathcal{S}\mbox{ and }\pimin_0\in\mathcal{S})=
\mathbf{P}(\pimin_0\in\mathcal{S})=
\mathsf{M}(\mathcal{S})\ge
\mathsf{M}(\mathcal{R}) \mathsf{M}_{\mathcal{R}}(\mathcal{S})>0.
\]
Next, we note that as $\mathbf{P}^i(Y_1=\tilde y_1,\ldots,Y_N=\tilde
y_N)>0$,
\[
\mathbf{P}^i\bigl((Y_1,\ldots,Y_N)\in\mathcal{O}\bigr) =
\delta_i M(\tilde y_1)\cdots M(\tilde y_N)1
\varphi(\{\tilde y_1\})\cdots\varphi(\{\tilde y_N\})>0,
\]
so $\varphi^{\otimes N}(\mathcal{O})>0$ and
$\mu M(\tilde y_1)\cdots M(\tilde y_N)1>0$ for all $\mu\in\mathcal{S}$
(this holds by the definition of $\mathcal{R}$ and as
$\mathcal{S}\subset\mathcal{R}$). Finally, by the triangle inequality,
\[
\sup_{\mu,\nu\in\mathcal{S}}
\biggl\|
\frac{\mu M(\tilde y_1)\cdots M(\tilde y_N)}{
\mu M(\tilde y_1)\cdots M(\tilde y_N)1}-
\frac{\nu M(\tilde y_1)\cdots M(\tilde y_N)}{
\nu M(\tilde y_1)\cdots M(\tilde y_N)1}
\biggr\|\le\varepsilon.
\]
Thus Condition \ref{condC} is satisfied, and the proof is complete.

\section{\texorpdfstring{Proof of Theorem \protect\ref{thm:mainf}}{Proof of Theorem 2.7}}
\label{sec:mainf}

The implication $3\Rightarrow1$ is already established by Theorem
\ref{thm:mainc}. It therefore suffices to prove the implications
$1\Rightarrow2$ and $2\Rightarrow3$.

\subsection{Proof of $1\Rightarrow2$}

We will need the following lemma.

\begin{lem}
\label{lem:singular}
Let $(\Xi,\mathcal{X},(\mathcal{X}_n)_{n\in\mathbb{N}},\boldsymbol{\Phi})$
be a filtered probability space, and let $\mathbf{Q},\mathbf{Q}'$ be
mutually singular probability measures on $\Xi$. Suppose that
$\mathbf{Q},\mathbf{Q}'$ are locally absolutely continuous with respect to
$\boldsymbol{\Phi}$, that is,
$\mathbf{Q}|_{\mathcal{X}_n}\ll\boldsymbol{\Phi}|_{\mathcal{X}_n}$ and
$\mathbf{Q}'|_{\mathcal{X}_n}\ll\boldsymbol{\Phi}|_{\mathcal{X}_n}$ with
densities $q_n$ and $q_n'$, respectively. Then
$q_n/q_n'\to0$ as $n\to\infty$ $\mathbf{Q}'$-a.s.
\end{lem}

\begin{pf}
Let $\boldsymbol{\tilde\Phi}=(\boldsymbol{\Phi}+\mathbf{Q}+\mathbf{Q}')/3$,
and let $r_n$ be the density of $\boldsymbol{\Phi}|_{\mathcal{X}_n}$
with respect to $\boldsymbol{\tilde\Phi}|_{\mathcal{X}_n}$.
Then we have $q_nr_n\to d\mathbf{Q}/d\boldsymbol{\tilde\Phi}$ and
$q_n'r_n\to d\mathbf{Q}'/d\boldsymbol{\tilde\Phi}$
$\boldsymbol{\tilde\Phi}$-a.s., hence $\mathbf{Q}'$-a.s.\ as
$\mathbf{Q}'\ll\boldsymbol{\tilde\Phi}$. But
$d\mathbf{Q}/d\boldsymbol{\tilde\Phi}=0$ $\mathbf{Q}'$-a.s.\ and
$d\mathbf{Q}'/d\boldsymbol{\tilde\Phi}>0$ $\mathbf{Q}'$-a.s.\
by the mutual singularity of $\mathbf{Q}$ and $\mathbf{Q}'$. The claim
follows directly.
\end{pf}

We consider the finite state space $E=\{1,\ldots,p\}$. Let us write
\[
\Omega_1,\ldots,\Omega_p\subseteq F^{\mathbb{N}},\qquad
\Omega_i = \supp\mathbf{P}^i|_{\mathcal{F}^Y_{1,\infty}}.
\]
There exists a finite partition $\{A_1,\ldots,A_K\}$ of $F^{\mathbb{N}}$
such that $\sigma\{A_1,\ldots,A_K\}=\sigma\{\Omega_1,\ldots,\Omega_p\}$.
We may assume without loss of generality that
\begin{eqnarray*}
\mathbf{P}^i\bigl((Y_k)_{k\in\mathbb{N}}\in A_1\bigr)&>&0\qquad
\mbox{for }i=1,\ldots,q,\\
\mathbf{P}^i\bigl((Y_k)_{k\in\mathbb{N}}\in A_1\bigr)&=&0\qquad
\mbox{for }i=q+1,\ldots,p
\end{eqnarray*}
for some $q\in E$ (this can always be accomplished by relabeling the
points of the state space). Define $\mathbf{\tilde P}( \cdot)=
\mathbf{P}^1( \cdot|(Y_k)_{k\in\mathbb{N}}\in A_1)$. Then, by
construction, $\mathbf{\tilde P}|_{\mathcal{F}^Y_{1,\infty}}\ll
\mathbf{P}^i|_{\mathcal{F}^Y_{1,\infty}}$ for $i=1,\ldots,q$ and
$\mathbf{\tilde P}|_{\mathcal{F}^Y_{1,\infty}}\perp
\mathbf{P}^i|_{\mathcal{F}^Y_{1,\infty}}$ for $i=q+1,\ldots,p$.

We assume that the filter is uniquely ergodic, so that
$\|\pi^x_n-\pi_n\|\to0$ $\mathbf{P}^x$-a.s.\ for every $x\in E$ by
Theorem \ref{thm:kunita} (this follows as $\lambda$ charges all points
in $E$, so we certainly have $\delta_x\ll\lambda$ for all $x\in E$).
Therefore, we find that
\[
\lim_{n\to\infty}\|\pi^x_n-\pi_n\|=0
\qquad \mbox{for }x=1,\ldots,q,\
\mathbf{\tilde P}\mbox{-a.s.}
\]
Denote by $q_n^i$ the density of
$\mathbf{P}^i((Y_1,\ldots,Y_n)\in\cdot)$ with respect to
$\varphi^{\otimes n}$, and similarly denote by $\tilde q_n$ the density of
$\mathbf{\tilde P}((Y_1,\ldots,Y_n)\in\cdot)$ with respect to
$\varphi^{\otimes n}$. Then
\[
\tilde q_n,q_n^x>0
\qquad \mbox{for all }n\in\mathbb{N}, x=1,\ldots,q
\quad \mbox{and}\quad
\lim_{n\to\infty}\frac{\tilde q_n}{q_n^1}<\infty
\qquad \mathbf{\tilde P}\mbox{-a.s.}
\]
(the latter follows as $\tilde q_n/q_n^1$ is a uniformly integrable
martingale under $\mathbf{P}^1$), while
\[
\lim_{n\to\infty}\frac{q_n^i}{\tilde q_n}=0
\qquad \mbox{for }i=q+1,\ldots,p,\
\mathbf{\tilde P}\mbox{-a.s.}
\]
by Lemma \ref{lem:singular}. (The fact that $\varphi$ may be any
$\sigma$-finite measure does not preclude the application of Lemma
\ref{lem:singular}, as $\varphi$ can always be transformed into a
probability measure by means of an equivalent change of measure.)

As all of the above statements hold $\mathbf{\tilde P}$-a.s., we can
certainly find one sample path $(\tilde y_k)_{k\in\mathbb{N}}$ on which
all these statements hold simultaneously. In particular, we have
\[
\biggl\|
\frac{\delta_xM(\tilde y_1)\cdots
M(\tilde y_n)}{\delta_xM(\tilde y_1)\cdots M(\tilde y_n)1}-
\frac{\lambda M(\tilde y_1)\cdots M(\tilde y_n)}{\lambda
M(\tilde y_1)\cdots M(\tilde y_n)1}
\biggr\|\xrightarrow{n\to\infty}0\qquad
\mbox{for }x=1,\ldots,q,
\]
as well as
\[
\frac{\delta_iM(\tilde y_1)\cdots M(\tilde y_n)1}{
\delta_1M(\tilde y_1)\cdots M(\tilde y_n)1} =
\frac{q_n^i(\tilde y_1,\ldots,\tilde y_n)}{
q_n^1(\tilde y_1,\ldots,\tilde y_n)}
\xrightarrow{n\to\infty}0\qquad
\mbox{for }i=q+1,\ldots,p.
\]
Now define the matrix norm $\|M\|:=\sup_{\|f\|_\infty\le1}
\sup_{\|\mu\|_1\le1}\mu Mf$. As the set of matrices of unit norm is
compact (as we are in the finite-dimensional setting), there must be a
subsequence $n_k\nearrow\infty$ and a matrix $M_\infty$ such that
\[
\frac{M(\tilde y_1)\cdots M(\tilde y_{n_k})}{
\|M(\tilde y_1)\cdots M(\tilde y_{n_k})\|}
\xrightarrow{k\to\infty}M_\infty,\qquad
\|M_\infty\|=1.
\]
We claim that $M_\infty$ is a rank 1 matrix. Indeed, for $i=q+1,\ldots,p$
we have
\[
\|\delta_iM_\infty\| =
\lim_{k\to\infty}\frac{\|\delta_iM(\tilde y_1)\cdots
M(\tilde y_{n_k})\|}{\|M(\tilde y_1)\cdots
M(\tilde y_{n_k})\|} \le
\lim_{k\to\infty}\frac{\delta_iM(\tilde y_1)\cdots
M(\tilde y_{n_k})1}{\delta_1M(\tilde y_1)\cdots
M(\tilde y_{n_k})1} = 0.
\]
On the other hand, consider a state $x\in\{1,\ldots,q\}$ such that
$\|\delta_xM_\infty\|>0$. Then $\delta_xM_\infty
1=\|\delta_xM_\infty\|>0$, and thus also $\lambda M_\infty1>0$.
But then
\begin{eqnarray*}
&&\biggl\|\frac{\delta_x M_\infty}{\delta_xM_\infty1}-
\frac{\lambda M_\infty}{\lambda M_\infty1} \biggr\|
\\
&&\qquad = \lim_{k\to\infty}
\biggl\|
\frac{\delta_xM(\tilde y_1)\cdots
M(\tilde y_{n_k})}{\delta_xM(\tilde y_1)\cdots M(\tilde
y_{n_k})1}-
\frac{\lambda M(\tilde y_1)\cdots M(\tilde y_{n_k})}{\lambda
M(\tilde y_1)\cdots M(\tilde y_{n_k})1} \biggr\|=0.
\end{eqnarray*}
Therefore, we have shown that for every $j=1,\ldots,p$, the $j$th row of
$M_\infty$ is either zero, or a multiple of the row vector $\lambda
M_\infty$. Moreover, $M_\infty$ is not identically zero as
$\|M_\infty\|=1$. Thus $M_\infty$ is a rank 1 matrix, and Condition \ref{condKR}
follows.

\subsection{Proof of $2\Rightarrow3$}

We assume that Condition \ref{condKR} holds. Therefore, there exists a
nonnegative column vector $u$ (which is not identically zero), and a
probability measure $\varrho$, such that the rank 1 matrix $u \varrho$ is
in the closure of the cone $\mathcal{C}$. In particular, for any
$\delta>0$, we can choose $N\in\mathbb{N}$, $y_1,\ldots,y_N\in F$,
$c>0$ such that
\[
\|c M(y_1)\cdots M(y_N)-u \varrho\|<\delta.
\]
Let $\alpha>0$ (to be chosen below), and define the set
\[
\mathcal{S}=\{\mu\in\mathcal{P}(E)\dvtx\mu u>\alpha\}.
\]
Then we can estimate
\[
\sup_{\mu\in\mathcal{S}}
\biggl\|\frac{\mu M(y_1)\cdots M(y_N)}{\mu M(y_1)\cdots M(y_N)1}
-\varrho
\biggr\| \le\frac{2\delta}{\alpha}.
\]
In particular, by the triangle inequality,
\[
\sup_{\mu,\nu\in\mathcal{S}}
\biggl\|
\frac{\mu M(y_1)\cdots M(y_N)}{\mu M(y_1)\cdots M(y_N)1}-
\frac{\nu M(y_1)\cdots M(y_N)}{\nu M(y_1)\cdots M(y_N)1}
\biggr\| \le\frac{4\delta}{\alpha}.
\]
Moreover, note that
\[
c \mu M(y_1)\cdots M(y_N)1 \ge
\mu u - \|c M(y_1)\cdots M(y_N)-u \varrho\|
> \alpha-\delta
\]
for all $\mu\in\mathcal{S}$.

We aim to show that Condition \ref{condC} is satisfied. To this end, let
$\varepsilon>0$ be given (and $\varepsilon<1$ without loss of generality).
As $\lambda u>0$, we may choose $\alpha=
\lambda u/2$ and $\delta=\alpha\varepsilon/4$. Now choose
$N\in\mathbb{N}$, $y_1,\ldots,y_N\in F$, $c>0$ as above. When $F$ is at
most countable, the above choices of $N$ and $\mathcal{S}$, together with
the singleton $\mathcal{O}=\{(y_1,\ldots,y_N)\}$, satisfy properties 1--3
of Condition \ref{condC}. Indeed, properties 2 and 3 are immediate from the above
computations. To prove property 1, note that
\begin{eqnarray*}
&&\mathbf{P}(\pimin_0\in\mathcal{S} \mbox{ and }
\pimax_0\in\mathcal{S})\\
&&\qquad \ge
\mathbf{P}(\pimax_0u\ge\pimin_0u>\lambda u/2) \\
&&\qquad =\mathbf{E}\bigl(\mathbf{P}(\pimax_0u\ge\pimin_0u|\mathcal{F}^Y_{-\infty,0})
I_{]\lambda u/2,\infty[}(\pimin_0u)\bigr).
\end{eqnarray*}
As trivially $\mathbf{P}(X\ge\mathbf{E}(X))>0$ for any random variable
$X$, we have $\mathbf{P}$-a.s.
\[
\mathbf{P}( \pimax_0u\ge\pimin_0u |\mathcal{F}^Y_{-\infty,0}) =
\mathbf{P}\bigl( \pimax_0u\ge
\mathbf{E}(\pimax_0u|\mathcal{F}^Y_{-\infty,0})
|\mathcal{F}^Y_{-\infty,0}\bigr)>0,
\]
while $\mathbf{P}(\pimin_0u>\lambda u/2)>0$ by virtue of the fact that
$\mathbf{E}(\pimin_0u)=\lambda u$. Therefore
$\mathbf{P}(\pimin_0\in\mathcal{S}\ \mathrm{and}\
\pimax_0\in\mathcal{S})>0$, and the claim is established.

In the case where $F=\mathbb{R}^d$, we cannot choose $\mathcal{O}$ to be
a singleton as this set has Lebesgue measure zero. However, note that by
the assumed continuity, all the above computations extend to a
sufficiently small neighborhood of the path $(y_1,\ldots,y_N)$.
Choosing $\mathcal{O}$ to be such a neighborhood, we have
$\varphi^{\otimes N}(\mathcal{O})>0$ by construction, and the remainder of
the proof proceeds as in the countable case.

\section{Sufficient conditions}
\label{sec:kaijser}

Our main results, Theorems \ref{thm:mainc} and \ref{thm:mainf}, establish
necessary and sufficient conditions for unique ergodicity of the filter.
The purpose of this section is to discuss various sufficient conditions
that have appeared in the literature, and their relations to our main
results. First, we discuss the sufficient conditions introduced by
Kaijser \cite{Kai75,Kai09} and show how these can be obtained directly
from our Theorem \ref{thm:mainc}. Then, we discuss various conditions
that have been introduced in the context of the filter stability problem
\cite{Van08,Van09,CV10}.

\subsection{Kaijser's sufficient conditions}

In Kaijser's 1975 paper \cite{Kai75}, the following condition is shown to
be sufficient for unique ergodicity of the filter.

\renewcommand{\thecondk}{(K)}
\begin{condk}\label{condK}
Let $E$ and $F$ be finite sets, let $\varphi$ be the counting measure on
$F$ and let the signal transition matrix $P$ be aperiodic. There exist
$y_1,\ldots,y_n\in F$ such that the matrix $M=M(y_1)\cdots M(y_n)$ is
nonzero and subrectangular, that is, $M_{ij}>0$ and $M_{kl}>0$ imply
$M_{il}>0$ and $M_{kj}>0$.
\end{condk}

Kaijser's proof of sufficiency is based on the Furstenberg--Kesten theory
of products of random matrices. A much simpler proof was given by Kochman
and Reeds in \cite{KR06}, Section 5, where Condition \ref{condK} is shown to imply
Condition \ref{condKR} through an explicit computation. Kochman and Reeds prove
the sufficiency of Condition \ref{condKR} by invoking a general result in the
ergodic theory of Markov chains in topological state spaces. We would
argue that the proof of sufficiency given here is even simpler, at least
if one takes for granted the (essentially known) characterization of
unique ergodicity of the filter provided by Theorem \ref{thm:kunita}.

Kaijser showed already in \cite{Kai75} by means of a counterexample that
the subrectangularity condition cannot be dropped, that is, that
irreducibility and aperiodicity of the signal need not imply unique
ergodicity of the filter. Kochman and Reeds provide two further
counterexamples \cite{KR06}. They demonstrate that the assumption of
aperiodicity cannot be dropped in Condition \ref{condK}, that is, that
subrectangularity and irreducibility need not imply unique ergodicity of
the filter. Moreover, they provide a counterexample where Condition \ref{condKR}
is satisfied and the signal is irreducible and aperiodic, but Condition \ref{condK} is not satisfied. Theorem \ref{thm:mainf} in this paper completes
these results by establishing the necessity of Condition \ref{condKR}.

In a recent paper, Kaijser \cite{Kai09} introduces two sufficient
conditions for unique ergodicity of the filter in the case where $E$ and
$F$ are countable.

\renewcommand{\thecondbb}{(B1)}
\begin{condbb}\label{condB1}
Let $E$ and $F$ be countable, and let $\varphi$ be the counting
measure.
There exists a nonnegative function $u\dvtx E\to\mathbb{R}_+$ with
$\|u\|_\infty=1$, a probability measure $\varrho$ on $E$, a sequence of
integers $(n_k)_{k\in\mathbb{N}}$ and a sequence of observation paths
$(y_1^k,\ldots,y_{n_k}^k)_{k\in\mathbb{N}}$ with $\|M(y_1^k)\cdots
M(y_{n_k}^k)\|>0$, such that
\[
 \biggl\|\frac{\delta_x M(y_1^k)\cdots M(y_{n_k}^k)}
 {\|M(y_1^k)\cdots M(y_{n_k}^k)\|}-u(x) \varrho\biggr\|
 \xrightarrow{k\to\infty}0\qquad
 \mbox{for all }x\in E.
\]
(Here we have defined the norm $\|M\|:=\sup_{\|f\|_\infty\le
1}\sup_{\|\mu\|_1\le1}\mu M f$.)
\end{condbb}

\renewcommand{\thecondb}{(B)}
\begin{condb}\label{condB}
Let $E$ and $F$ be countable, and let $\varphi$ be the counting
measure. For every $\beta>0$, there exists an $x_0\in E$ such that
the following holds: given any tight set
$\mathcal{T}\subset\mathcal{P}(E)$ such that, for any
$\mathsf{M}_0\in\mathcal{P}(\mathcal{P}(E))$ with
$\int\nu\mathsf{M}_0(d\nu) = \lambda$,
\[
\mathsf{M}_0\bigl(\mathcal{T}\cap\{\nu\in\mathcal{P}(E)\dvtx
\nu_{x_0}>\lambda_{x_0}/2\}\bigr) \ge\lambda_{x_0}/3,
\]
there exist $N\in\mathbb{N}$ and $y_1,\ldots,y_N\in F$ such that
$\delta_{x_0}M(y_1)\cdots M(y_N)1>0$ and
\[
\biggl\|
\frac{\mu M(y_1)\cdots M(y_N)}{\mu M(y_1)\cdots M(y_N)1} -
\frac{\delta_{x_0} M(y_1)\cdots M(y_N)}{\delta_{x_0}
M(y_1)\cdots M(y_N)1}
\biggr\|<\beta
\]
for all $\mu\in\mathcal{T}\cap\{\nu\in\mathcal{P}(E)\dvtx
\nu_{x_0}>\lambda_{x_0}/2\}$.
\end{condb}

Kaijser shows that either of these conditions implies unique ergodicity of
the filter, provided the signal transition matrix $P$ is aperiodic.
Kaijser's proof is very long and requires the development of some
dedicated ergodicity results for Markov chains in nonlocally compact
spaces. We will presently show that Condition \ref{condB1} and Condition \ref{condB}
imply our Condition \ref{condC}, so that Kaijser's results follow easily from
Theorem \ref{thm:mainc} (even in the case where $P$ is not aperiodic).

\begin{lem}
\label{lem:condb1}
Condition \ref{condB1} implies Condition \ref{condC}.
\end{lem}

\begin{pf}
Suppose that Condition \ref{condB1} holds.
We can estimate
\begin{eqnarray*}
&&\biggl\|\frac{\mu M(y_1^k)\cdots M(y_{n_k}^k)}
{\|M(y_1^k)\cdots M(y_{n_k}^k)\|}-\mu u \varrho\biggr\|\\
&&\qquad \le
\int\biggl\|\frac{\delta_x M(y_1^k)\cdots M(y_{n_k}^k)}
{\|M(y_1^k)\cdots M(y_{n_k}^k)\|}-u(x) \varrho\biggr\|
\mu(dx) \\
&&\qquad \le \sum_{x=1}^J \mu_x
\biggl\|\frac{\delta_x M(y_1^k)\cdots M(y_{n_k}^k)}
{\|M(y_1^k)\cdots M(y_{n_k}^k)\|}-u(x) \varrho\biggr\|
+2\sum_{x=J+1}^\infty\mu_x.
\end{eqnarray*}
Let $\mathcal{T}\subset\mathcal{P}(E)$ be a tight set. Then the first
term converges to zero uniformly in $\mu\in\mathcal{T}$ by assumption,
while the second term can be made arbitrarily small uniformly in
$\mu\in\mathcal{T}$ by choosing $J$ sufficiently large. Therefore,
\[
\sup_{\mu\in\mathcal{T}}
\biggl\|\frac{\mu M(y_1^k)\cdots M(y_{n_k}^k)}
{\|M(y_1^k)\cdots M(y_{n_k}^k)\|}-\mu u \varrho\biggr\|
<\delta
\]
for any tight set $\mathcal{T}\subset\mathcal{P}(E)$, $\delta>0$, and
$k$ sufficiently large. Let $\alpha>0$ and define
\[
\mathcal{S}=\mathcal{T}\cap\{\mu\in\mathcal{P}(E)\dvtx\mu u>\alpha\}.
\]
Then we obtain
\[
\sup_{\mu,\nu\in\mathcal{S}}
\biggl\|\frac{\mu M(y_1^k)\cdots M(y_{n_k}^k)}
{\mu M(y_1^k)\cdots M(y_{n_k}^k)1}-
\frac{\nu M(y_1^k)\cdots M(y_{n_k}^k)}
{\nu M(y_1^k)\cdots M(y_{n_k}^k)1}
\biggr\|
\le\frac{4\delta}{\alpha}.
\]
We now show that Condition \ref{condC} is satisfied. Let $\varepsilon>0$ be
given, and choose $\alpha= \lambda u/2$ and $\delta
=\alpha\varepsilon/4$. As in the proof of Theorem \ref{thm:mainf}, we
can show that
\[
\mathbf{P}(\pimin_0u>\alpha\mbox{ and }\pimax_0u>\alpha)>0.
\]
Moreover, we can find an
increasing sequence of tight sets $\mathcal{T}_n\subset\mathcal{P}(E)$
such that
\[
\mathbf{P}(\pimin_0\in\mathcal{T}_n\mbox{ and }
\pimax_0\in\mathcal{T}_n)\xrightarrow{n\to\infty}1,
\]
as $\mathcal{P}(\mathcal{P}(E))$ is Polish.
Therefore, we can choose $\mathcal{T}$ sufficiently large such that
\[
\mathbf{P}(\pimin_0\in\mathcal{S}\mbox{ and }
\pimax_0\in\mathcal{S})>0.
\]
The remainder of the proof is identical to that of
Theorem \ref{thm:mainf}.
\end{pf}

\begin{lem}
Condition \ref{condB} implies Condition \ref{condC}.
\end{lem}

\begin{pf}
Suppose that Condition \ref{condB} holds. We claim that Condition \ref{condC} holds with
$\varepsilon= 2\beta$, $\mathcal{S}=
\mathcal{T}\cap\{\nu\in\mathcal{P}(E)\dvtx\nu_{x_0}>\lambda_{x_0}/2\}$,
and $\mathcal{O}=\{(y_1,\ldots,y_N)\}$, provided that
$\mathcal{T}\subset\mathcal{P}(E)$ is chosen sufficiently large.
Indeed, as the family
$\mathcal{M}=\{\mathsf{M}_0\in\mathcal{P}(\mathcal{P}(E))\dvtx
\int\nu \mathsf{M}_0(d\nu)=\lambda\}$ is tight (e.g.,
\cite{Jak88}) it is easily seen that
\[
\mathsf{M}_0\bigl(\mathcal{T}\cap\{\nu\in\mathcal{P}(E)\dvtx
\nu_{x_0}>\lambda_{x_0}/2\}\bigr) \ge\lambda_{x_0}/3\qquad
\mbox{for all }\mathsf{M}_0\in\mathcal{M}
\]
is satisfied for every sufficiently large tight set
$\mathcal{T}\subset\mathcal{P}(E)$. Moreover,
\[
\mathbf{P}(\pimin_0\in\mathcal{S}\mbox{ and }
\pimax_0\in\mathcal{S})>0
\]
when $\mathcal{T}$ is chosen sufficiently large,
as is shown in the proof of Lemma \ref{lem:condb1}.
It remains to note that as $\delta_{x_0}M(y_1)\cdots M(y_N)1>0$,
we have $\mu M(y_1)\cdots M(y_N)1>0$ for all $\mu\in\mathcal{S}$.
The remainder of Condition \ref{condC} now follows immediately.
\end{pf}

Though Condition \ref{condB1} is strongly reminiscent of Condition \ref{condKR}, we did
not succeed in extending the proof of the necessity of Condition \ref{condKR}
to the
countable case. Whether Conditions \ref{condB1}, \ref{condB} or some variant of
thereof are necessary and sufficient for unique ergodicity in the
countable case remains an open problem.

\subsection{Nondegeneracy and observability}

Conditions of a rather different kind than are considered by Kaijser,
Kochman and Reeds relate to the filter stability problem
(see the survey \cite{CV10}). By Theorem \ref{thm:kunita}, however, filter
stability and unique ergodicity are essentially equivalent, so that also
these conditions can be brought to bear on the problem considered in this
paper. In this section, we consider the following conditions that are
borrowed
from \cite{Van08,Van09}: nondegeneracy [Condition~\ref{condN}], uniform
observability [Condition \ref{condUO}] and observability [Condition~\ref{condO}].

\renewcommand{\thecondn}{(N)}
\begin{condn}\label{condN}
If $i,j\in E$ and $P_{ij}>0$, then $M_{ij}(y)>0$ for all $y\in F$.
\end{condn}

\renewcommand{\theconduo}{(UO)}
\begin{conduo}\label{condUO}
For every $\varepsilon>0$, there is a $\delta>0$ such that
\[
\|\mathbf{P}^\mu|_{\mathcal{F}^Y_{1,\infty}}-
\mathbf{P}^\nu|_{\mathcal{F}^Y_{1,\infty}}\|<\delta
\quad \mbox{implies}\quad
\|\mu-\nu\|<\varepsilon\qquad
\mbox{[for any }\mu,\nu\in\mathcal{P}(E)\mbox{]}.
\]
\end{conduo}

\renewcommand{\thecondo}{(O)}
\begin{condo}\label{condO}
If $\mu,\nu\in\mathcal{P}(E)$ and
$\mathbf{P}^\mu|_{\mathcal{F}^Y_{1,\infty}}=
\mathbf{P}^\nu|_{\mathcal{F}^Y_{1,\infty}}$, then $\mu=\nu$.
\end{condo}

\begin{thm}
\label{thm:maino}
Suppose that one of the following holds:
\begin{itemize}
\item Condition \ref{condN} holds, and the signal transition matrix $P$ is
aperiodic; or
\item Condition \ref{condUO} holds; or
\item Condition \ref{condO} holds, and $E$ is a finite set.
\end{itemize}
Then the filter admits a unique invariant measure $\mathsf{M}$, and
$n^{-1}\sum_{k=1}^n\mathsf{M}_0\mathsf{\Pi}^k\Rightarrow\mathsf{M}$
as $n\to\infty$ for any $\mathsf{M}_0\in\mathcal{P}(\mathcal{P}(E))$.
If, in addition, the signal transition matrix $P$ is aperiodic, then
we have $\mathsf{M}_0\mathsf{\Pi}^n\Rightarrow\mathsf{M}$ as $n\to\infty$
for any $\mathsf{M}_0\in\mathcal{P}(\mathcal{P}(E))$.
\end{thm}

\begin{pf*}{Sketch of proof}
First, suppose that Condition \ref{condN} holds and that $P$ is aperiodic.
Consider the stochastic process
$(\mathbf{X}_n,\mathbf{Y}_n)_{n\in\mathbb{Z}}$ defined as
\[
\mathbf{X}_n = (X_{n},X_{n+1})\in\mathbf{E},\qquad
\mathbf{Y}_n = Y_{n+1}
\in F,
\]
where $\mathbf{E}=\{x\in E^2\dvtx\mathbf{P}(\mathbf{X}_0=x)>0\}$.
Then $(\mathbf{X}_n,\mathbf{Y}_n)_{n\in\mathbb{Z}}$ is a stationary Markov
chain, $(\mathbf{X}_n)_{n\in\mathbb{Z}}$ is an irreducible and aperiodic
Markov chain, and $(\mathbf{Y}_n)_{n\in\mathbb{Z}}$ are conditionally
independent given $(\mathbf{X}_n)_{n\in\mathbb{Z}}$.
Moreover,
\[
\mathbf{P}\bigl(\mathbf{Y}_n\in A|
(\mathbf{X}_k)_{k\in\mathbb{Z}}\bigr)
=
\int_A M_{X_{n}X_{n+1}}(y)
\varphi(dy)
:=
\int_A\Upsilon(\mathbf{X}_n,y) \varphi(dy),
\]
where $\Upsilon(x,y)>0$ for all $x\in\mathbf{E}$ and $y\in F$
by Condition \ref{condN}. Therefore,
\[
\|\mathbf{P}^\mu(\mathbf{X}_n\in\cdot|
\mathbf{Y}_0,\ldots,\mathbf{Y}_n)-
\mathbf{P}^\nu(\mathbf{X}_n\in\cdot|
\mathbf{Y}_0,\ldots,\mathbf{Y}_n)\|
\xrightarrow{n\to\infty}0\qquad
\mathbf{P}^\mu\mbox{-a.s.}
\]
for all $\mu,\nu\in\mathcal{P}(E)$ by \cite{Van08}, Corollary 5.5. It
follows immediately that $\|\pi^\mu_{n}-\pi^\nu_{n}\|\to0$ as
$n\to\infty$ $\mathbf{P}^\mu$-a.s. The proof is completed by invoking
Theorem \ref{thm:kunita}.

Next, suppose Condition \ref{condUO} holds. By a result of Blackwell and
Dubins \cite{BD62},
\[
\bigl\|\mathbf{P}^\mu\bigl((Y_k)_{k>n}\in\cdot|\mathcal{F}_{1,n}^Y\bigr)-
\mathbf{P}^\nu\bigl((Y_k)_{k>n}\in\cdot|\mathcal{F}_{1,n}^Y\bigr)\bigr\|
\xrightarrow{n\to\infty}0\qquad
\mathbf{P}^\mu\mbox{-a.s.}
\]
whenever $\mu\ll\nu$. But one can show (e.g., \cite{Van09}) that
\[
\mathbf{P}^\rho\bigl((Y_k)_{k>n}\in\cdot|\mathcal{F}_{1,n}^Y\bigr) =
\mathbf{P}^{\pi_n^\rho}\bigl((Y_k)_{k>0}\in\cdot\bigr) =
\mathbf{P}^{\pi_n^\rho}|_{\mathcal{F}_{1,\infty}^Y}\qquad
\mbox{for all }\rho\in\mathcal{P}(E).
\]
Using Condition \ref{condUO}, it therefore follows that
$\|\pi_n^\mu-\pi_n^\nu\|\to0$ as $n\to\infty$ $\mathbf{P}^\mu$-a.s.\
whenever $\mu\ll\nu$. The proof is completed by invoking Theorem
\ref{thm:kunita}.

Finally, suppose that $E$ is finite, and that Condition \ref{condO} holds. Then it
is not difficult to establish, along the lines of \cite{Van09},
Proposition 3.5, that Condition \ref{condUO} is satisfied. The result therefore
follows as above.
\end{pf*}

\begin{rem}
When the signal transition kernel $P$ is periodic, Condition \ref{condN} by itself
does not ensure unique ergodicity of the filter (this can be seen,
e.g., by considering the example of a periodic signal in $E=\{1,2\}$
with the trivial observation state space $F=\{1\}$). However, if
Condition \ref{condN} holds and $E$ is a finite set, a \emph{detectability}
condition [which is weaker than Condition \ref{condO}] is necessary and sufficient
for stability of the filter, and hence for unique ergodicity. The
necessary arguments can be adapted from \cite{Van09ptrf}, Section 6.2, with
some care. As this is somewhat outside the scope of this paper, we omit
the details.
\end{rem}

It should be noted that none of the conditions of Theorem \ref{thm:maino}
are necessary. Indeed, Condition \ref{condN} is not satisfied by the examples
given by Kochman and Reeds \cite{KR06}. That Condition \ref{condUO} [hence
Condition \ref{condO}] is not necessary can be seen from the trivial
counterexample, where $P$ is aperiodic and $F=\{1\}$. In this case
the observations are completely noninformative, so that the point mass at
$\lambda\in\mathcal{P}(E)$ is the unique invariant measure for the filter,
but Condition \ref{condUO} is not satisfied.

Nonetheless, the sufficient conditions of Theorem \ref{thm:maino} can be
useful in practice, as they may be substantially easier to check than
Condition \ref{condC} or \ref{condKR}. For example, in the case where $E$ is a finite
set, verifying Condition \ref{condO} is simply a matter of linear algebra (see
\cite{CV10} for an example), while verifying Condition \ref{condKR} involves
taking limits. Moreover, despite that Conditions \ref{condC} and \ref{condKR} are both
necessary and sufficient in many cases, we did not succeed in our attempt
to prove Theorem \ref{thm:maino} by directly verifying that Condition \ref{condC} or \ref{condKR}
hold. Therefore, such sufficient but not necessary conditions
remain of independent interest.

\begin{appendix}

\section*{Appendix: Supplementary proofs}
\label{app:proofs}

\setcounter{thm}{0}

\subsection{\texorpdfstring{Proof of Lemmas \protect\ref{lem:filtrecurs} and
\protect\ref{lem:filtmarkov}}{Proof of Lemmas 2.3 and 2.4}}
\label{sec:basicfiltering}

We will need the following.

\begin{lem}
\label{lem:shiftminmax}
$\pimin_n=\pimin_m\circ\Theta^{n-m}$
and $\pimax_n=\pimax_m\circ\Theta^{n-m}$
for $m,n\in\mathbb{Z}$.
\end{lem}

\begin{pf}
By stationarity of $\mathbf{P}$, it is easily seen that
\begin{eqnarray*}
\mathbf{E}(f(X_n)|\mathcal{F}_{n-\ell,n}^Y) &=&
\mathbf{E}(f(X_m)|\mathcal{F}_{m-\ell,m}^Y)\circ\Theta^{n-m},
\\
\mathbf{E}\bigl(f(X_n)|\mathcal{F}_{n-\ell,n}^Y\vee
\mathcal{F}_{n-\ell,n-k}^X\bigr) &=&
\mathbf{E}\bigl(f(X_m)|\mathcal{F}_{m-\ell,m}^Y\vee
\mathcal{F}_{m-\ell,m-k}^X\bigr)\circ\Theta^{n-m}.
\end{eqnarray*}
The result follows by letting $\ell\to\infty$, then $k\to\infty$.
\end{pf}

We begin by proving Lemma \ref{lem:filtrecurs} for the
case $\pi_n^\mu$. It clearly suffices to prove
\[
\pi_n^\mu=
\frac{\mu M(Y_1)\cdots M(Y_n)}{
\mu M(Y_1)\cdots M(Y_n)1},\qquad
\mathbf{P}^\mu\mbox{-a.s. for all }n\ge1.
\]
Let $f\in\mathcal{C}_b(E)$ and $A\in\mathcal{B}(F^n)$. Then
\begin{eqnarray*}
&&\mathbf{E}^\mu\biggl(
\frac{\mu M(Y_1)\cdots M(Y_n)f}{
\mu M(Y_1)\cdots M(Y_n)1} I_A(Y_1,\ldots,Y_n)
\biggr) \\
&&\qquad =
\int_A
\frac{\mu M(y_1)\cdots M(y_n)f}{
\mu M(y_1)\cdots M(y_n)1}
\mu M(y_1)\cdots M(y_n)1 \varphi(dy_1)\cdots\varphi(dy_n) \\
&&\qquad =
\int_A
\mu M(y_1)\cdots M(y_n)f \varphi(dy_1)\cdots\varphi(dy_n) \\
&&\qquad =
\mathbf{E}^\mu(I_A(Y_1,\ldots,Y_n) f(X_n)).
\end{eqnarray*}
As this holds for any $f\in\mathcal{C}_b(E)$ and $A\in\mathcal{B}(F^n)$,
the above expression for $\pi_n^\mu$ follows from the definition of the
conditional expectation.

To prove Lemma \ref{lem:filtrecurs} for $\pimin_n$, let $k,n\ge1$.
Note that
\begin{eqnarray*}
\mathbf{E}(f(X_n)|\mathcal{F}^Y_{-k+1,n}) &=&
\pi_{n+k}f \circ\Theta^{-k} \\
&=&
\frac{\lambda M(Y_{-k+1})\cdots M(Y_n)f}{
\lambda M(Y_{-k+1})\cdots M(Y_n)1} \\
&=&
\frac{(\pi_k\circ\Theta^{-k}) M(Y_1)\cdots M(Y_n)f}{
(\pi_k\circ\Theta^{-k}) M(Y_1)\cdots M(Y_n)1}.
\end{eqnarray*}
But $\mathbf{E}(f(X_n)|\mathcal{F}^Y_{-k+1,n})\to
\pimin_nf$ and $\pi_kf\circ\Theta^{-k}=
\mathbf{E}(f(X_0)|\mathcal{F}^Y_{-k+1,0})\to\pimin_0f$ as $k\to\infty$
$\mathbf{P}$-a.s.\ by the martingale convergence theorem. Therefore
\[
\pimin_n =
\frac{\pimin_0 M(Y_1)\cdots M(Y_n)}{
\pimin_0 M(Y_1)\cdots M(Y_n)1},\qquad
\mathbf{P}\mbox{-a.s.} \mbox{ for all }n\ge1,
\]
and the result follows for arbitrary $m,n\in\mathbb{Z}$, $n\ge m$ by Lemma
\ref{lem:shiftminmax}.

To prove Lemma \ref{lem:filtrecurs} for $\pimax_n$, let $n\ge1$
and $k\ge\ell\ge0$. Note that
\begin{eqnarray*}
\mathbf{E}\bigl(f(X_n)|\mathcal{F}^Y_{-k,n}\vee
\mathcal{F}^X_{-k,-\ell}\bigr)& =&
\mathbf{E}\bigl(f(X_{n+\ell})|\mathcal{F}^Y_{-k+\ell,n+\ell}\vee
\mathcal{F}^X_{-k+\ell,0}\bigr)\circ\Theta^{-\ell}
\\
& =&
\mathbf{E}\bigl(f(X_{n+\ell})|\mathcal{F}^Y_{1,n+\ell}\vee
\sigma\{X_0\}\bigr)\circ\Theta^{-\ell},
\end{eqnarray*}
where we have used the Markov property. Moreover, it is easily seen that
\[
\mathbf{E}\bigl(f(X_{n+\ell})|\mathcal{F}^Y_{1,n+\ell}\vee
\sigma\{X_0\}\bigr) =
\pi_{n+\ell}^{X_0}f =
\frac{\delta_{X_0} M(Y_1)\cdots M(Y_{n+\ell})f}{
\delta_{X_0} M(Y_1)\cdots M(Y_{n+\ell})1}.
\]
Therefore, we can write
\begin{eqnarray*}
\mathbf{E}\bigl(f(X_n)|\mathcal{F}^Y_{-k,n}\vee
\mathcal{F}^X_{-k,-\ell}\bigr)
& =&
\frac{\delta_{X_{-\ell}} M(Y_{-\ell+1})\cdots M(Y_n)f}{
\delta_{X_{-\ell}} M(Y_{-\ell+1})\cdots M(Y_n)1} \\
&=&
\frac{(\pi^{X_0}_\ell\circ\Theta^{-\ell}) M(Y_{1})\cdots M(Y_n)f}{
(\pi^{X_0}_\ell\circ\Theta^{-\ell})M(Y_{1})\cdots M(Y_n)1}.
\end{eqnarray*}
Letting $k\to\infty$, then $\ell\to\infty$ and applying the martingale
convergence theorem, we obtain the desired recursion for $\pimax_n$.

We now turn to the proof of Lemma \ref{lem:filtmarkov}. The stationarity
of $(\pimax_n)_{n\in\mathbb{Z}}$ and $(\pimin_n)_{n\in\mathbb{Z}}$ follows
directly from Lemma \ref{lem:shiftminmax} and the stationarity of
$(X_n,Y_n)_{n\in\mathbb{Z}}$. It only remains to prove the Markov
property. For $f\in\mathcal{C}_b(\mathcal{P}(E))$, we can compute
\begin{eqnarray*}
\mathbf{E}(f(\pimax_{n+1})|\mathcal{G}_{-\infty,n}) &=&
\mathbf{E} \biggl( f \biggl(\frac{\pimax_{n}M(Y_{n+1})}{
\pimax_n M(Y_{n+1})1} \biggr) \bigg|\mathcal{G}_{-\infty,n} \biggr)
\\
&=&
\mathbf{E} \biggl( f \biggl(\frac{\mu M(Y_{n+1})}{
\mu M(Y_{n+1})1} \biggr) \bigg|\mathcal{G}_{-\infty,n} \biggr)
\bigg|_{\mu=\pimax_n},
\end{eqnarray*}
where we have used Lemma \ref{lem:filtrecurs} and the fact that
$\pimax_n$ is $\mathcal{G}_{-\infty,n}$-measurable. But for any
bounded measurable function $g\dvtx F\to\mathbb{R}$, we have
\[
\mathbf{E}(g(Y_{n+1})|\mathcal{G}_{-\infty,n}) =
\int g(y) \pimax_nM(y)1 \varphi(dy).
\]
The Markov property and the expression for the transition kernel
$\mathsf{\Pi}$ follows immediately. The Markov property of
$\pimin_n$ and $\pi^\mu_n$ follows along similar lines.

\subsection{\texorpdfstring{Proof of Theorem \protect\ref{thm:kunita} and Lemma
\protect\ref{lem:limminmax}}{Proof of Theorem 3.1 and Lemma 3.4}}
\label{sec:kunita}

The proof of Theorem \ref{thm:kunita} follows closely along the lines of
\cite{Kun71,Ste89,DS05}. We will sketch the necessary arguments,
concentrating on the special features of the countable setting.

We begin by establishing the Feller property.

\begin{lem}
Let $(\mu_n)_{n\in\mathbb{N}}\subset\mathcal{P}(E)$ and
$\mu\in\mathcal{P}(E)$ be such that $\mu_n\Rightarrow\mu$. Then $\int
f(\nu) \mathsf{\Pi}(\mu_n,d\nu)\to\int f(\nu) \mathsf{\Pi}(\mu,d\nu)$
for every $f\in\mathcal{C}_b(\mathcal{P}(E))$.
\end{lem}

\begin{pf}
Let $N\subset F$ be a $\varphi$-null set such that $\sup_i\sum_j
M_{ij}(y)<\infty$ for all $y\notin N$. Then $\mu_nM(y)1\to\mu M(y)1$
for all $y\notin N$, and $\mu_nM(y)/\mu_nM(y)1\Rightarrow\mu M(y)/\mu
M(y)1$ whenever $y\notin N$ and $\mu M(y)1>0$. It follows that
\[
f \biggl(\frac{\mu_n M(y)}{\mu_n M(y)1} \biggr)
\mu_n M(y)1 \xrightarrow{n\to\infty}
f \biggl(\frac{\mu M(y)}{\mu M(y)1} \biggr)
\mu M(y)1 \qquad \mbox{for all }y\notin N.
\]
But the family $\{f(\mu_n M(y)/\mu_n M(y)1) \mu_n M(y)1\dvtx n\in\mathbb{N}\}$
is uniformly integrable (under $\varphi$), as
$|f(\mu_n M(y)/\mu_n M(y)1) \mu_n M(y)1|\le\|f\|_\infty
\mu_n M(y)1$ and by Scheff{\'e}'s lemma $\int|\mu_n M(y)1-\mu
M(y)1| \varphi(dy)\to0$. The result therefore follows from the
expression for $\mathsf{\Pi}$ in Lemma \ref{lem:filtmarkov}.
\end{pf}

We will need some basic elements from Choquet theory.

\begin{defn}
Let $S$ be Polish. For
$\mathsf{M},\mathsf{M}'\in\mathcal{P}(\mathcal{P}(S))$
we write $\mathsf{M}\prec\mathsf{M}'$ if
\[
\int f(\nu) \mathsf{M}(d\nu) \le
\int f(\nu) \mathsf{M}'(d\nu)\qquad
\mbox{for every convex }
f\in\mathcal{C}_b(\mathcal{P}(S)).
\]
For any $\mathsf{M}\in\mathcal{P}(\mathcal{P}(S))$, the
\emph{barycenter} $b(\mathsf{M})\in\mathcal{P}(S)$ is defined as
\[
b(\mathsf{M})u = \int\nu u \mathsf{M}(d\nu)\qquad
\mbox{for all }u\in\mathcal{C}_b(S).
\]
For any $\mu\in\mathcal{P}(S)$, define
$\mathsf{m}_\mu,\mathsf{\tilde m}_\mu\in\mathcal{P}(\mathcal{P}(S))$ as
\[
\int f(\nu) \mathsf{m}_\mu(d\nu) =
f(\mu),\qquad
\int f(\nu) \mathsf{\tilde m}_\mu(d\nu) =
\int f(\delta_x) \mu(dx)
\]
for every $f\in\mathcal{C}_b(\mathcal{P}(S))$.
\end{defn}

\begin{lem}
\label{lem:choq}
Let $S$ be a Polish space. The following hold:
\begin{enumerate}[(3)]
\item[(1)] Given $\mathsf{M}\in\mathcal{P}(\mathcal{P}(S))$, we have
$f(b(\mathsf{M}))\le\int f(\nu) \mathsf{M}(d\nu)$
for every convex function $f\in\mathcal{C}_b(\mathcal{P}(S))$
(Jensen's inequality).
\item[(2)] For any
$\mathsf{M}\in\mathcal{P}(\mathcal{P}(S))$, we have
$\mathsf{m}_{b(\mathsf{M})} \prec\mathsf{M} \prec
\mathsf{\tilde m}_{b(\mathsf{M})}$.
\item[(3)]
If $\mathsf{M},\mathsf{M}'\in\mathcal{P}(\mathcal{P}(S))$,
$\mathsf{M}\prec\mathsf{M}'$ and $\mathsf{M}'\prec\mathsf{M}$,
we have $\mathsf{M}=\mathsf{M}'$.
\end{enumerate}
In particular, $\prec$ defines a partial order on
$\mathcal{P}(\mathcal{P}(S))$.
\end{lem}

\begin{pf}
Jensen's inequality is proved as in \cite{Kun71}, Lemma 3.1.
The second property follows easily from Jensen's inequality.
The third property follows from the fact that the family of convex
functions in $\mathcal{C}_b(\mathcal{P}(S))$ is a measure determining
class (see, e.g., Proposition A1 in \cite{Ste89}).
\end{pf}

We now need some basic convexity properties of the filter.

\begin{lem}
\label{lem:convex}
The following hold for any $\mathsf{M}\in\mathcal{P}(\mathcal{P}(E))$:
\begin{enumerate}[(3)]
\item[(1)] If $f\in\mathcal{C}_b(\mathcal{P}(E))$ is convex, then
$\mathsf{\Pi}f\in\mathcal{C}_b(\mathcal{P}(E))$ is also convex.
\item[(2)]$b(\mathsf{M}\mathsf{\Pi})=b(\mathsf{M})P$.
\item[(3)] If $\mathsf{M}\mathsf{\Pi}=\mathsf{M}$, then
$b(\mathsf{M})=\lambda$.
\item[(4)]$\mathsf{m}_{b(\mathsf{M})P^n}\mathsf{\Pi}^{m}
\prec\mathsf{M}\mathsf{\Pi}^{m+n}
\prec\mathsf{\tilde m}_{b(\mathsf{M})P^n}\mathsf{\Pi}^{m}$
for any $m,n\ge0$.
\end{enumerate}
\end{lem}

\begin{pf}
The first claim follows as in \cite{Kun71}, Lemma 3.2. The second claim
follows directly from Lemma \ref{lem:filtmarkov}. The third claim
follows from the second claim and the fact that $\lambda$ is the unique
invariant measure for $P$. The fourth claim follows from the first and
second claims, together with the second claim of Lemma \ref{lem:choq}.
\end{pf}

The following lemma connects $\pimin_0,\pimax_0$ to the filter transition
kernel $\mathsf{\Pi}$.

\begin{lem}
\label{lem:weakminmax}
Denote by $\mathsf{M}^{\max},\mathsf{M}^{\min}\in\mathcal{P}(\mathcal{P}(E))$ the laws of
$\pimax_0$ and $\pimin_0$, respectively. Then
$\mathsf{m}_\lambda\mathsf{\Pi}^n\Rightarrow\mathsf{M}^{\min}$
and $\mathsf{\tilde m}_\lambda\mathsf{\Pi}^n\Rightarrow
\mathsf{M}^{\max}$ as $n\to\infty$.
\end{lem}

\begin{pf}
Let $f\in\mathcal{C}_b(\mathcal{P}(E))$. Then
\[
\mathsf{m}_\lambda\mathsf{\Pi}^nf =
\mathbf{E}[f(\pi_n)] =
\mathbf{E}\bigl[f\bigl(\mathbf{P}(X_0\in\cdot|\mathcal{F}^Y_{-n+1,0})\bigr)\bigr]
\xrightarrow{n\to\infty}\mathbf{E}[f(\pimin_0)],
\]
where we have used stationarity and martingale convergence. Similarly,
\begin{eqnarray*}
\mathsf{\tilde m}_\lambda\mathsf{\Pi}^nf &=&
\mathbf{E}[f(\pi_n^{X_0})] =
\mathbf{E}\bigl[f\bigl(\mathbf{P}(X_0\in\cdot|\mathcal{F}^Y_{-n+1,0}
\vee\sigma\{X_{-n}\})\bigr)\bigr]  \\
&=&
\mathbf{E}\bigl[f\bigl(\mathbf{P}(X_0\in\cdot|\mathcal{F}^Y_{-\infty,0}
\vee\mathcal{F}^X_{-\infty,-n})\bigr)\bigr]
\xrightarrow{n\to\infty}\mathbf{E}[f(\pimax_0)],
\end{eqnarray*}
where we have additionally used the Markov property of
$(X_n,Y_n)_{n\in\mathbb{Z}}$.
\end{pf}

Finally, we will need the following convergence property:

\begin{lem}
\label{lem:totvar}
$\lim_{n\to\infty}\|\pi_n^\mu-\pi_n^\nu\|$ exists
$\mathbf{P}^\mu$-a.s.\ whenever $\mu\ll\nu$.
\end{lem}

\begin{pf}
It is not difficult to show along the lines of \cite{Van08}, Corollary
5.7, that
\begin{eqnarray*}
&&\|\pi_n^\mu-\pi_n^\nu\|
\\
&&\qquad =
\frac{\mathbf{E}^\nu( |\mathbf{E}^\nu(
(d\mu/d\nu)(X_0)|\mathcal{F}^X_{n,\infty}\vee
\mathcal{F}^Y_{1,\infty}) - \mathbf{E}^\nu(
(d\mu/d\nu)(X_0)|\mathcal{F}^Y_{1,n})|
|\mathcal{F}^Y_{1,n})}{
\mathbf{E}^\nu((d\mu/d\nu)(X_0)|\mathcal{F}^Y_{1,n})
}
\end{eqnarray*}
$\mathbf{P}^\mu$-a.s.\ whenever $\mu\ll\nu$. The denominator converges
$\mathbf{P}^\nu$-a.s., hence $\mathbf{P}^\mu$-a.s.\ (as $\mu\ll\nu$), to
a random variable which is strictly positive $\mathbf{P}^\mu$-a.s.

To prove convergence of the numerator, let $\varepsilon>0$ and define
\begin{eqnarray*}
M_n &=& \mathbf{E}^\nu\biggl(\frac{d\mu}{d\nu}(X_0)
I_{d\mu/d\nu(X_0)<\varepsilon}
\big|\mathcal{F}^Y_{1,n}\biggr),\\
M_n' &=& \mathbf{E}^\nu\biggl(\frac{d\mu}{d\nu}(X_0)
I_{d\mu/d\nu(X_0)\ge\varepsilon}
\big|\mathcal{F}^Y_{1,n}\biggr),\\
L_n &=& \mathbf{E}^\nu\biggl(\frac{d\mu}{d\nu}(X_0)
I_{d\mu/d\nu(X_0)<\varepsilon}
\big|\mathcal{F}^X_{n,\infty}\vee
\mathcal{F}^Y_{1,\infty}\biggr),\\
L_n' &=& \mathbf{E}^\nu\biggl(\frac{d\mu}{d\nu}(X_0)
I_{d\mu/d\nu(X_0)\ge\varepsilon}
\big|\mathcal{F}^X_{n,\infty}\vee
\mathcal{F}^Y_{1,\infty}\biggr).
\end{eqnarray*}
Clearly $M_n$ and $M_n'$ are uniformly integrable martingales, while
$L_n$ and $L_n'$ are reverse martingales. Moreover, the numerator
can be written as $\mathbf{E}^\nu(Z_n|\mathcal{F}^Y_{1,n})$ where
$Z_n=|L_n+L_n'-M_n-M_n'|$. We proceed to estimate as follows:
\[
|\mathbf{E}^\nu(Z_n-Z_\infty|\mathcal{F}^Y_{1,n})| \le
\mathbf{E}^\nu(|L_n-L_\infty| |\mathcal{F}^Y_{1,n}) +
\mathbf{E}^\nu(|M_n-M_\infty| |\mathcal{F}^Y_{1,n}) +
4M_n'.
\]
The first two terms converge to zero $\mathbf{P}^\nu$-a.s.\ as
$n\to\infty$ by Hunt's lemma (\cite{BD62}, Theorem 2), while
$\lim_{n\to\infty}M_n'$ vanishes if we let $\varepsilon\to\infty$.
Therefore $\mathbf{E}^\nu(Z_n-Z_\infty|\mathcal{F}^Y_{1,n})\to0$ as
$n\to\infty$ $\mathbf{P}^\nu$-a.s., and the proof is easily completed.
\end{pf}

We now proceed to the proof of Theorem \ref{thm:kunita} and
Lemma \ref{lem:limminmax}.

\subsubsection{\texorpdfstring{Proof of Theorem \protect\ref{thm:kunita}
$(1\Leftrightarrow2)$}{Proof of Theorem 3.1 $(1\Leftrightarrow2)$}}

First suppose $\mathbf{P}(\pimax_0=\pimin_0)<1$. Then
$\mathbf{E}(\{(\pimax_0)_i-(\pimin_0)_i\}^2)>0$ for some $i\in E$.
Now note that
\begin{eqnarray*}
\mathbf{E}\bigl(\{(\pimax_0)_i-(\pimin_0)_i\}^2\bigr) &=&
\mathbf{E}(\{(\pimax_0)_i\}^2) -
\mathbf{E}(\{(\pimin_0)_i\}^2) \\
& =&\int\{\nu_i\}^2 \mathsf{M}^{\max}(d\nu) -
\int\{\nu_i\}^2 \mathsf{M}^{\min}(d\nu),
\end{eqnarray*}
so that $\mathbf{P}(\pimax_0=\pimin_0)<1$ implies $\mathsf{M}^{\max}\ne
\mathsf{M}^{\min}$. But $\mathsf{M}^{\max}$ and $\mathsf{M}^{\min}$ are invariant measures for $\mathsf{\Pi}$ by Lemma
\ref{lem:filtmarkov}, so we have shown that the filter admits two distinct
invariant measures. Conversely, if the invariant measure is unique,
then $\mathbf{P}(\pimax_0=\pimin_0)=1$. Thus we have proved the
implication $1\Rightarrow2$.

Now suppose that $\pimax_0=\pimin_0$, so that in particular
$\mathsf{M}^{\max}=\mathsf{M}^{\min}$. Let $\mathsf{M}$ be any
invariant measure for $\mathsf{\Pi}$. We claim that $\mathsf{M}^{\min}\prec\mathsf{M}\prec\mathsf{M}^{\max}$, so that necessarily
$\mathsf{M}=\mathsf{M}^{\max}=\mathsf{M}^{\min}$ by Lemma
\ref{lem:choq}. To prove the claim, note that
$\mathsf{m}_\lambda\mathsf{\Pi}^n\prec\mathsf{M}\mathsf{\Pi}^n=
\mathsf{M}\prec\mathsf{\tilde m}_\lambda\mathsf{\Pi}^n$ for any $n\ge0$
by Lemmas \ref{lem:choq} and \ref{lem:convex}. The claim therefore
follows directly from Lemma \ref{lem:weakminmax}. Thus we have proved the
implication $2\Rightarrow1$.

\subsubsection{\texorpdfstring{Proof of Theorem \protect\ref{thm:kunita}
$(2\Leftrightarrow3)$}{Proof of Theorem 3.1
$(2\Leftrightarrow3)$}}

Proceeding along the same lines as in the proof of Lemma
\ref{lem:weakminmax} (and taking into account the fact that weak
convergence and total variation convergence of probability measures
coincide when the state space is countable), one can show that
\[
\mathbf{E}(\|\pi_n^{X_0}-\pi_n\|)
\xrightarrow{n\to\infty}\mathbf{E}(\|\pimax_0-\pimin_0\|).
\]
Suppose first that property 3 holds. Then
\[
\mathbf{E}(\|\pi_n^{X_0}-\pi_n\|) =
\sum_{i\in E}\lambda_i \mathbf{E}^i(\|\pi_n^i-\pi_n\|)
\xrightarrow{n\to\infty}0.
\]
Therefore $\pimax_0=\pimin_0$, and we have proved the
implication $3\Rightarrow2$.

Conversely, suppose that $\pimax_0=\pimin_0$. Let
$\mu,\nu\in\mathcal{P}(E)$ such that $\mu\ll\nu$.
Note that we can write
$\pi^\mu_n=\mathbf{E}^\mu(\pi^{X_0}_n|\mathcal{F}^Y_{1,n})$.
Therefore, we have
\[
\mathbf{E}^\mu(\|\pi^\mu_n-\pi_n\|) =
\mathbf{E}^\mu\bigl(\|\mathbf{E}^\mu(\pi^{X_0}_n-\pi_n|
\mathcal{F}^Y_{1,n})\|\bigr) \le
\mathbf{E}^\mu(\|\pi^{X_0}_n-\pi_n\|).
\]
But $\mathbf{E}(\|\pi^{X_0}_n-\pi_n\|)\to0$, so $\|\pi^{X_0}_n-\pi_n\|
\to
0$ in probability. As $\mu\ll\lambda$, we find that
$\|\pi^{X_0}_n-\pi_n\|\to0$ in $\mathbf{P}^\mu$-probability also, and by
dominated convergence
\[
\mathbf{E}^\mu(\|\pi^\mu_n-\pi_n\|) \le
\mathbf{E}^\mu(\|\pi^{X_0}_n-\pi_n\|)
\xrightarrow{n\to\infty}0.
\]
By Lemma \ref{lem:totvar}, it follows that $\|\pi^\mu_n-\pi_n\|\to0$
$\mathbf{P}^\mu$-a.s. Similarly, we find that $\|\pi^\nu_n-\pi_n\|\to0$
$\mathbf{P}^\nu$-a.s., hence $\mathbf{P}^\mu$-a.s.\ as
$\mu\ll\nu$. Therefore $\|\pi_n^\mu-\pi_n^\nu\|\to0$
$\mathbf{P}^\mu$-a.s., and we have evidently proved the implication
$2\Rightarrow3$.

\subsubsection{\texorpdfstring{Proof of Theorem \protect\ref{thm:kunita}
$(1\Leftrightarrow4)$}{Proof of Theorem 3.1
$(1\Leftrightarrow4)$}}

The implication $4\Rightarrow1$ follows immediately by choosing
$\mathsf{M}_0$ and $\mathsf{M}$ to be distinct invariant measures of the
filter and applying property 4, which leads to a contradiction.

To prove the converse implication, choose $\mathsf{M}_0$ arbitrarily and
define the measures
$\mathsf{M}_n=n^{-1}\sum_{k=1}^n\mathsf{M}_0\mathsf{\Pi}^n$. Note that
$b(\mathsf{M}_n)=n^{-1}\sum_{k=1}^nb(\mathsf{M}_0)P^n \Rightarrow\lambda$
as the signal is irreducible and positive recurrent. It follows from
\cite{Jak88} that the sequence $(\mathsf{M}_n)_{n\in\mathbb{N}}$ is tight.
It therefore suffices to prove that every convergent subsequence has the
same limit. But it is easily seen that any convergent subsequence
converges to an invariant measure of $\mathsf{\Pi}$, so that the result
follows from the uniqueness of the invariant measure. Thus we have proved
the implication $1\Rightarrow4$.

\subsubsection{\texorpdfstring{Proof of Theorem \protect\ref{thm:kunita}
$(1\Leftrightarrow5)$}{Proof of Theorem 3.1
$(1\Leftrightarrow5)$}}

The implication $5\Rightarrow1$ follows immediately by choosing
$\mathsf{M}_0$ and $\mathsf{M}$ to be distinct invariant measures of the
filter and applying property 5, which leads to a contradiction.

We prove the converse implication under the assumption that the signal
transition matrix $P$ is aperiodic. Choose $\mathsf{M}_0$ arbitrarily,
and note that $b(\mathsf{M}_0\mathsf{\Pi}^n)=b(\mathsf{M}_0)P^n
\Rightarrow\lambda$. It follows from \cite{Jak88} that the sequence
$(\mathsf{M}_0\mathsf{\Pi}^n)_{n\in\mathbb{N}}$ is tight. It therefore
suffices to prove that any convergent subsequence converges to the
unique invariant measure of the filter $\mathsf{M}$. Let
$n(k)\nearrow\infty$ be a subsequence such that
$\mathsf{M}_0\mathsf{\Pi}^{n(k)}\Rightarrow\mathsf{M}_\infty$, and let
$f\in\mathcal{C}_b(\mathcal{P}(E))$ be convex. By Lemma \ref{lem:convex},
we have
\[
\mathsf{m}_{b(\mathsf{M}_0)P^{n(k)-m}}\mathsf{\Pi}^{m}f
\le\mathsf{M}_0\mathsf{\Pi}^{n(k)}f
\le\mathsf{\tilde m}_{b(\mathsf{M}_0)P^{n(k)-m}}\mathsf{\Pi}^{m}f\qquad
\mbox{for all }m\le n(k).
\]
In particular, letting $k\to\infty$ and using the Feller property, we have
\[
\mathsf{m}_{\lambda}\mathsf{\Pi}^{m}f
\le\mathsf{M}_\infty f
\le\mathsf{\tilde m}_{\lambda}\mathsf{\Pi}^{m}f\qquad
\mbox{for all }m\ge0.
\]
But letting $m\to\infty$ and using Lemma \ref{lem:weakminmax}, we find
that $\mathsf{M}^{\min}\prec\mathsf{M}_\infty\prec\mathsf{M}^{\max}$. As the invariant measure $\mathsf{M}$ is presumed to be unique, we
have $\mathsf{M}^{\min}=\mathsf{M}^{\max}=\mathsf{M}$ by the
implication $1\Rightarrow2$. Therefore, we find that
$\mathsf{M}_\infty=\mathsf{M}$ by Lemma \ref{lem:choq}.
This completes the proof of the implication $1\Rightarrow5$.

\subsubsection{\texorpdfstring{Proof of Lemma
\protect\ref{lem:limminmax}}{Proof of Lemma 3.4}}

First, we note that $\mathbf{E}((\pimax_0)_i|\mathcal{F}^Y_{-\infty,0})=
(\pimin_0)_i$. Therefore, $\mathbf{P}((\pimin_0)_i=0\mbox{ and }
(\pimax_0)_i>0)=0$ for every $i\in E$. In particular, this implies that
we have $\pimax_0\ll\pimin_0$ with unit probability under $\mathbf{P}$.
Now note that $\pimax_k=\pi_k^\mu|_{\mu=\pimax_0}$ and
$\pimin_k=\pi_k^\mu|_{\mu=\pimin_0}$ by Lemma \ref{lem:filtrecurs}.
Therefore,
\begin{eqnarray*}
&&\mathbf{P} \Bigl( \lim_{k\to\infty} \|\pimax_k
-\pimin_k \|\mbox{ exists} \big|
\mathcal{G}_{-\infty,0} \Bigr) \\
&&\qquad =
\mathbf{P} \Bigl( \lim_{k\to\infty} \|\pi_k^\mu
-\pi_k^\nu\|\mbox{ exists} \big|
\mathcal{G}_{-\infty,0} \Bigr)
\Big|_{\mu=\pimax_0, \nu=\pimin_0} \\
&&\qquad =
\mathbf{P}^\mu\Bigl(\lim_{k\to\infty} \|\pi_k^\mu
-\pi_k^\nu\|\mbox{ exists} \Bigr)
\Big|_{\mu=\pimax_0, \nu=\pimin_0} =
1 \qquad \mathbf{P}\mbox{-a.s.},
\end{eqnarray*}
where we have used the fact that $\pimin_0$ and $\pimax_0$ are
$\mathcal{G}_{-\infty,0}$-measurable in the first step, the Markov
property in the second step, and Lemma \ref{lem:totvar} in the third step.
The result now follows by taking the expectation with respect to
$\mathbf{P}$.
\end{appendix}

%

\printaddresses

\end{document}